\documentclass[notitlepage,leqno,11pt]{article}

\usepackage{amsfonts}
\usepackage{latexsym}
\usepackage{amssymb}
\usepackage{amsmath}

\newtheorem{Theorem}{Theorem}[section]
\newtheorem{Lemma}[Theorem]{Lemma}
\newtheorem{Proposition}[Theorem]{Proposition}
\newtheorem{Corollary}[Theorem]{Corollary}
\newtheorem{Remark}[Theorem]{Remark}
\newtheorem{Definition}[Theorem]{Definition}

\def\Proof{\noindent{{\it Proof. }}}
\def\square{\vbox{
    \hrule height .4pt
    \hbox{\vrule width .4pt height 7pt \kern 7pt
       \vrule width .4pt}
    \hrule height .4pt }}
\def\QED{\hfill {$\square$}\goodbreak \medskip}

\def\R{{\mathbb R}}
\def\RN{{\R^n}}
\def\N{{\mathbb N}}

\newcommand{\var }{\varphi }

\def\e{{\varepsilon}}
\def\a{{\alpha}} 
\def\b{{\beta}}  
\def\d{{\delta}} 
\def\l{{\lambda}} 
\def\n{{\nabla}}  
\def\r{{\rho}}  
\def\s{{\sigma}}  
\def\g{{\gamma}} 
\def\o{{\omega}} 

\newcommand{\ov}{\overline}
\title{\Large Adiabatic limits of closed orbits for some \\ 
Newtonian systems in $\R^n$}

\author{Andrea Malchiodi\footnote{Supported by M.U.R.S.T., Variational 
Methods and Nonlinear Differential Equations}}
\date{}

\begin{document}

\smallskip

\maketitle

\centerline{\small{
Scuola Internazionale Superiore di Studi Avanzati}}
\centerline{\small{via Beirut, 2-4 -- 34014 Trieste, Italy}}

\bigskip\bigskip

{\footnotesize
\begin{abstract}

\noindent
We deal with a Newtonian system like $\ddot{x} + V'(x) = 0$. 	
We suppose that $V: \R^n \to \R$ possesses an 
$(n-1)$-dimensional compact manifold $M$ of critical points, and we prove
the existence of arbitrarity slow periodic orbits. When the period tends 
to infinity these orbits, rescaled in time, converge 
to some closed geodesics on $M$.    
 
\bigskip\bigskip

	
\end{abstract}
}

\small
\section{Introduction}

Let $V : \R^n \to \R$ be a smooth function and suppose that $V$ possesses an 
$(n-1)$-dimensional compact manifold $M$ of critical points which is non
degenerate, namely 
\begin{equation}\label{eq:ndv} 
ker V''(x) = T_xM \quad \forall x \in M, \qquad 
\mbox{ or equivalently } \qquad V''[n_x,n_x] \neq 0 \quad \forall x \in M,
\end{equation}
where $n_x$ is a normal vector to $M$ at $x$. 

We are interested in studying the existence of solutions to the Newtonian
system  
\begin{equation}\tag{$T$}\label{eq:per-T} 
\begin{cases} 
\ddot{x} + V'(x) = 0; &\\
x(\cdot) \mbox{ is } T-\mbox{periodic}, &
\end{cases}
\end{equation} 
when $T$ is large and $x(\cdot)$ is close to $M$.  
Equivalently, setting $\varepsilon^2 = \frac{1}{T}$, one looks for solutions
to  the problem  
 
\begin{equation}\tag{$P_\e$}\label{eq:pe} 
\begin{cases} 
\ddot{x} + \frac{1}{\e} V'(x) = 0;& \\
x(\cdot) \mbox{ is } 1-\mbox{periodic},&
\end{cases}
\end{equation} 
for $\varepsilon >0 $ sufficiently small. 

\ 

\noindent
As $T$ varies, problem (\ref{eq:per-T}) can possess 
some continuous families of solutions parametrized in $T$, and the fact that
$V$ is degenerate (in the sense of Morse) allows these solutions to have a 
non trivial limit behaviour. The case of solutions approaching a critical 
manifold of $V$ has been considered for example in \cite{k}, \cite{m} and
\cite{m2}. It is known that if some smooth family $x(t,\e)$ solves 
problem (\ref{eq:pe}) and if $x(t,\e) \to M$ as $\e \to 0$, then the curve 
$x(t,0)$ is a geodesic on $M$. The curve $x(\cdot,0)$ is called adiabatic
limit  for the family $x(\cdot,\e)$. 

The aim of this paper is to achieve some complemetary result, namely to prove 
that for some closed geodesics on $M$ there are indeed solutions of
(\ref{eq:pe}) which approach these geodesics. Large
period orbits with some limit behaviour have also been studied in \cite{ci}
for planar systems.   

Our main results are the following Theorems. The first one treats the case
of a non-degenerate closed geodesic on $M$, see Definition \ref{d:nd}.

\begin{Theorem}\label{t:non-deg}
Suppose $x_0 : S^1 \to M$ is a non-degenerate closed geodesic, and suppose
$V$ is repulsive w.r.t. $M$, namely that the following condition holds
\begin{equation}\label{eq:v-rep}
V''(x)[n_x,n_x] < 0 \qquad \mbox{ for all } x \in M. 
\end{equation}
Then there exists $T_0 > 0$ with the following property. For all $T \geq T_0$
there exists a function $u_T$ such that 
\begin{description}
\item{(i)} $u_T$ is solution of problem (\ref{eq:per-T});  
\item{(ii)} as $T \to +\infty$, $u_T (T \cdot) \to x_0(\cdot )$ in
$C^1(S^1;\R^n)$.     
\end{description}
\end{Theorem}
The proof relies on the Local Inversion Theorem, which can be applied by the
non-degeneracy of $x_0$. Since in \eqref{eq:pe} appears the
singular term  $\frac{1}{\e}$, a quite accurate expansion of $V'$ is
needed.

If we want to prove the convergenge of a sequence of trajectories, instead of
the convergence of a one-parameter family, then we can remove any
non-degeneracy assumption. With abuse of notation we will again call
adiabatic limit the limit trajectory.

\begin{Theorem}\label{t:es-deg} 
If condition (\ref{eq:v-rep}) holds, then for every sequence $T_k \to +\infty$
there exists a sequence of solutions $(u_k)_k$ to problem \eqref{eq:per-T}
corresponding to $T = T_k$ such that up to subsequence $u_k(T_k \cdot)$
converge in $C^0(S^1,\R^n)$. The adiabatic limit of $u_k(T_k \cdot)$ is a 
non trivial closed geodesic $x_0$ on $M$.      
\end{Theorem} 

\begin{Remark}\label{r:i}
(a) Since the adiabatic limit $x_0$ in Theorem \ref{t:es-deg} can be
degenerate, and so it is possible that it belongs to a family of degenerate
geodesics, it is natural to expect convergence only along sequences of
trajectories. 

(b) The limit geodesic $x_0$ can be characterized as follows. If $\pi_1 (M)
\neq 0$, $x_0$ realizes the infimum of the square lenght in some component of
the closed loops in $M$. From the proof of Theorem \ref{t:es-deg} it follows
that one can find adiabatic limits belonging to every element of $\pi_1 (M)$,
see also Remark \ref{r:cl}. If $\pi_1 (M) = 0$ then the energy of $x_0$ is
the infimum in some suitable min-max scheme. 
\end{Remark}
The proof of Theorem \ref{t:es-deg} is based on a Lyapunov-Schmidt reduction
on the  Hilbert manifold $H^1(S^1; M)$ of the closed loops in $M$ of class
$H^1$. Standard min-max arguments are applied to a suitable functional
on $H^1(S^1; M)$  which is a perturbation of the square lenght $L_0$, 
see formula \eqref{eq:l0}. A similar approach has been used for example in 
\cite{ab} to perform reductions on finite-dimensional manifolds. The new 
feature of our method is that we perform a reduction of an infinite dimensional 
manifold.  

\

\noindent
If $V$ is of attractive type, namely if $V''(x)[n_x,n_x] > 0$ for all $x \in
M$,  then the situation is very different, since some phenomena of resonance 
may occur. As a consequence our hypotheses become stronger and we can 
prove convergence just for some suitable sequence $T_k \to +\infty$. Section
6 contains some results concerning this case. As an example we can state
the following one.

\begin{Theorem}\label{t:attr-es}
Suppose $V$ satisfies the following conditions for some $b_0 > 0$
\begin{description}
\item{(i)} $V''(x)[n_x,n_x] = b_0$ \quad for all  $x \in M$; 
\item{(ii)} $\frac{\partial^3 V}{\partial n_x^3}(x) = 0$ \quad
for all  $x \in M$.  
\end{description}
Then there exists a sequence $T_k \to +\infty$ and 
there exists a sequence of solutions $(u_k)_k$ to problem \eqref{eq:per-T}
corresponding to $T = T_k$ such that up to subsequence $u_k(T_k \cdot)$
converge in $C^0(S^1,\R^n)$. The adiabatic limit of $u_k(T_k \cdot)$ is a 
non trivial closed geodesic $x_0$ on $M$.    
\end{Theorem}
The paper is organized as follows: Section 2 is devoted to recalling some 
notations and preliminary facts. In Section 3 we prove Theorem
\ref{t:non-deg}. In Section 4 we study some linear ordinary differential 
equations, used to perform the reduction. In Section 5 we reduce the problem
on   $H^1(S^1; M)$, we study the reduced functional and we prove Theorem 
\ref{t:es-deg}. Finally in Section 6 we treat the attractive case.

\begin{center}
{\bf Acknowledgements}
\end{center}

The author whishes to thank Professor A. Ambrosetti for having proposed him 
the study of this problem and for his useful advices.

\section{Notations and Preliminaries}

In this Section we recall some well known facts in Riemannian
Geometry, we refer to \cite{k} or \cite{s} for the details. In particular we
introduce the Levi-Civita connection, the Gauss' equations, the Hilbert
manifold $H^1(S^1; M)$, and some properties of the square lenght functional
$L_0$ on $H^1(S^1; M)$.

It is given an orientable manifold $M \subseteq \R^n$ of codimension 1, which 
inherits naturally a Riemannian structure from $\R^n$. On $M$ it is defined
the Gauss map $n: M \to S^{n-1}$ which assigns to every point $x \in M$ the
unit versor $n_x \in (T_x M)^\perp$, where $T_x M$ is the tangent space of
$M$ at $x$. The differential of $n_x$ is given by 
\begin{equation}\label{eq:dn}
d \, n_x [v] = H(x) [v]; \qquad \forall x \in M, \quad \forall  v \in T_x M,
\end{equation}
where $H(x) : T_x M \to T_x M$ is a symmetric operator. Dealing with the
operator $H(x)$ we will also identify it with the corresponding symmetric
bilinear form according to the relation 
$$
H(x)[v,w] = H(x) [v] \cdot w ; \qquad \forall v, w \in T_x M. 
$$
The bilinear form $H(x)$ is called the {\em second fundamental form}  of $M$
at $x$.     

If $\mathcal{X}(M)$ denotes the class of the smooth vector fields on
the manifold $M$, then for every $x \in M$ the {\em Levi-Civita
connection } $\n : T_x M \times \mathcal{X}(M) \to T_x M$ is defined in the
following way 
\begin{equation}\label{eq:l-c}
\n_X Y = D_X \tilde{Y} - (n_x \cdot D_X \tilde{Y})\, n_x. 
\end{equation}
Here $\tilde{Y}$ is an extension of $Y$ in a neighbourhood of $x$ in $\R^n$
and $D_X$ denotes the standard differentiation in $\R^n$ along the 
direction $X$. 

\begin{Remark}\label{r:c-der}
The quantity $\n_X Y$ defined in \eqref{eq:l-c} is nothing but the projection
of $D_X \tilde{Y}$ onto the tangent space $T_x M$. Actually $\n_X Y$
depends only on $Y(x)$ and on the derivative of $Y$ along the direction $X$.
Hence formula (\ref{eq:l-c}) makes sense  also when $Y$ is defined just on a
curve $c$ on $M$ for which $c(t_0) = x$ and   $\dot{c}(t_0) = X$. In the
following this fact will be considered understood.   
\end{Remark}
The {\em Riemann tensor} $R : T_x M \times T_x M \times T_x M \to T_x M$ is
defined by  $$
\nabla_{\tilde{X}} \nabla_{\tilde{Y}} \tilde{Z} - 
\nabla_{\tilde{Y}} \nabla_{\tilde{X}} \tilde{Z} - 
\nabla_{[\tilde{X},\tilde{Y}]} \tilde{Z}. 
$$
Here $\tilde{X}, \tilde{Y}$ and $\tilde{Z}$ are smooth extensions of $X, Y$ and $Z$ 
respectively, and the symbol $[\cdot, \cdot]$ denotes the usual Lie bracket. 
The above definition does not depend on the extensions of $X, Y$ and $Z$.

Let $i : M \to \R^n$ be the inclusion of $M$ in $\R^n$: there exists a 
symmetric bilinear form $s : T_x M \times T_x M \to \R$ which satisfies  
$$ 
s_x(X,Y) = \perp \, D_{X} \tilde{Y}, \qquad \forall  X, Y \in T_x M.   
$$
Here $\tilde{Y}$ is any smooth extension of $Y$ and $\perp \n'_{X} \tilde{Y}$
is the component of $D_{X} \tilde{Y}$ normal to $T_x M$.  
The {\em Gauss' equations} are the following   
\begin{equation}\label{eq:for-g-f}
(R(X,Y)Z,W) = (s_x(Y,Z), s_x(X,W)) - (s_x(X,Z), s_x(Y,W)), \qquad \forall 
X, Y, Z, W \in T_x M.   
\end{equation}  
Given $X, Y \in T_x M$, consider a smooth extension $\tilde{Y}$ of $Y$ and an  
extension $\Psi$ of $n_x$ such that $\Psi(p) \perp T_p M$ for all $p
\in M$ in  some neighbourhood of $x$. Differentiating the relation
$(Y,\Psi) = 0$ along the  direction $X$, there results 
$$
0 = D_X (Y,\Psi) = (D_X Y, \Psi) + (Y, D_X \Psi). 
$$
Hence it follows that $s_x(X, Y)$ is given by
$$
s_x(X, Y) = - H(x)[X,Y], \qquad \forall X, Y \in T_x M.    
$$
In particular equation \eqref{eq:for-g-f} becomes
\begin{equation}\label{eq:for-g-f2}
(R(X,Y)Z,W) = (H(x)(Y,Z), H(x)(X,W)) - (H(x)(X,Z), H(x)(Y,W)), \qquad \forall 
X, Y, Z, W \in T_x M.   
\end{equation}  
In order to define the manifold $H^1(S^1;M)$, we recall first the
differentiable structure of a smooth $k$-dimensional manifold. This is given
by a family of local charts $(U_\a,\var_\a)_\a$, where  $\var_\a : \R^k \to
\R$  are diffeomorphisms such that the compositions 
$\var_\b^{-1} \circ \var_\a : \var_\a^{-1}(U_\a) \cap \var_\b^{-1}(U_\b)
\subseteq \R^k \to \R^k$ are smooth functions. 

\begin{Definition}
A closed curve $c : S^1 \to M$ is said to be of class $H^1(S^1;M)$ if for some 
chart $(U,\var)$ of $M$ the map $\var \circ c : S^1 \to \R^k$ is of class
$H^1$.   
\end{Definition} 
This definition does not depend on the choice of the
chart $(U,\var)$, since the composition of an $H^1$ map in $\R^k$ with a
smooth  diffeomorphism is still of class $H^1$. The class $H^1(S^1;M)$
constitutes an infinite dimensional Hilbert manifold. We recall briefly its
structure.  Given a curve $c \in C^\infty(S^1; M)$ we denote by $c^* TM$ the 
pull-back of $TM$ trough $c$, namely the family of vector fields $X :
S^1 \to TM$ such that  
$$
X(t) \in T_{c(t)} M\qquad \mbox{ for all } t \in S^1. 
$$ 
We define also $\mathcal{H}_c$ to be the sections $\xi$ of $c^* TM$ for which   
$$
\int_{S^1} |\xi(t)|^2 dt + \int_{S^1} |\n_{\dot{c}(t)} \xi(t)|^2 < + \infty. 
$$
There is a neighbourhood $\mathcal{U}$ of the zero section of $TM$ where the
exponential  map $exp : TM \to M$ is well defined. For $\xi \in \mathcal{U}
\cap \mathcal{H}_c$,  the curve $exp \, \xi$ belongs to $H^1(S^1; M)$, and
viceversa every curve in  $H^1(S^1; M)$ can be obtained in this way for a
suitable $c \in C^\infty(S^1; M)$.  Hence the family $(\mathcal{U} \cap
\mathcal{H}_c, exp)$, $c \in C^\infty(S^1; M)$  constitutes an atlas for
$H^1(S^1; M)$. 

The tangent space of $H^1(S^1; M)$ can be described as follows:  
if $h \in H^1(S^1; M)$, we consider the class of curves $\xi(\cdot)$ such that
$\xi(t) \in T_{h(t)} M$ for all $t \in S^1$ and such that  
$$
\int_{S^1} |\xi(t)|^2 dt + \int_{S^1} |\n_{\dot{h}(t)} \xi(t)|^2 < + \infty. 
$$
By means of the H\"older inequality one can define the scalar product on $T_h
H^1(S^1; M)$ as  
$$
(\xi,\eta)_1 = \int_{S^1} (\xi(t),\eta(t)) dt + 
\int_{S^1} (\nabla_{\dot{h}(t)} \xi(t),\nabla_{\dot{h}(t)} \eta(t)) dt,
\qquad \forall  \xi, \eta \in T_h H^1(S^1; M).
$$
This scalar product determines a positive definite bilinear form on
$T_h H^1(S^1; M)$ and hence a Riemannian structure on $H^1(S^1; M)$.

On the manifold $H^1(S^1; M)$ is defined the square lenght functional $L_0$
in the following way 
\begin{equation}\label{eq:l0} 
L_0(u) = \frac{1}{2}\int_0^1 |\dot{u}(t)|^2 dt, \qquad   
\mbox{ for all } u \in H^1(S^1; M). 
\end{equation}
Given $A > 0$ we define $L_0^A \subseteq H^1(S^1; M)$ to be 
$$
L_0^A = \{ u \in H^1(S^1; M) : L_0(u) \leq A \}.
$$
It is a standard fact that the functional $L_0$ is smooth on the
manifold $H^1(S^1; M)$ endowed with the above structure, and there results 
\begin{equation}\label{eq:der-l0} 
D L_0(h)[k] = \int_{S^1} (\dot{h},
\nabla_{\dot{h}} k); \qquad h \in H^1(S^1; M),  k \in T_h H^1(S^1; M).
\end{equation}
The critical points of $L_0$ are precisely the closed geodescics on $M$.  
Furthermore, if $h \in H^1(S^1; M)$ is a stationary point of $L_0$ there
results 
\begin{equation}\label{eq:de2r-l0}
D^2 L_0(h)[k,w] = \int_{S^1} (\nabla_{\dot{h}} k, \nabla_{\dot{h}}w) -  
\int_{S^1} R(k, \dot{x}_0, \dot{x}_0, w); \qquad k, w \in T_h H^1(S^1; M),
\end{equation}
where $R$ is given by formula \eqref{eq:for-g-f2}. 

\begin{Definition}\label{d:nd}
A closed geodesic $h \in H^1(S^1; M)$ is said to be non-degenerate if the
kernel of $D^2 L_0(h)$ is one dimensional, and hence coincides with the 
span of $\dot{h} \in T_h H^1(S^1; M)$.  
\end{Definition}
One useful property of $L_0$ is the following 
\begin{equation}\tag{$PS$}\label{eq:ps}
\mbox{ the functional } L_0 \mbox{ satisfies the Palais Smale condition on } 
H^1(S^1; M),
\end{equation}
namely every sequence $(u_m)$ for which $L_0(u_m) \to a \in \R$ and 
$\|DL_0(u_n)\| \to 0$ admits a convergent subsequence.

Condition (\ref{eq:ps}) allows to apply the standard min-max arguments in order
to prove existence of critical point of $L_0$. For example, if for some
compact manifold $M$ it is $\pi_1 (M) \neq 0$, we can reason as follows. 
The negative gradient flow of
$L_0$ preserves the components of $H^1(S^1;M)$ and the embedding $H^1(S^1; M)
\hookrightarrow C^0(S^1; M)$ is continuous. As an immediate aplication we
have the following Theorem. 

\begin{Theorem}
Let $M$ be a compact manifold. Then for every $\alpha \in \pi_1(M)$ there 
exists a non-constant closed geodesic $u_\alpha:[0,1]\to M$ such that 
$\left[ u_\alpha(\cdot)\right] = \alpha$, and moreover 
$$ 
E(u_\alpha) = c_\alpha := \inf_{\left[ u \right] = \alpha} E(u). 
$$ 
\end{Theorem}
If $M$ is simply connected, the proof of the existence of a closed geodesic is
more involved and in its most general form it is due to Lusternik and Fet, see
\cite{lf}, by means of topological methods. The proof of our Theorem
\ref{t:es-deg} follows that argument, and we will recall it later. A
fundamental tool is the Hurewitz Theorem.    

\begin{Theorem}\label{t:hur} (Hurewitz) 
Suppose $M$ is a finite dimensional compact manifold such that $\pi_1(M) = 0$. 
Define $q$ to be the smallest integer for which $\pi_q(M) \neq 0$, and define 
$q'$ to be the smallest integer such that $H_{q'}(M) \neq 0$. Then $q$ and 
$q'$ are equal. 
\end{Theorem}

\begin{Remark}\label{r:hur}
In our case $M$ is orientable, and there always results $H_{n-1}(M) \simeq 
\mathbb{Z}$, so it turns out that $q \leq n - 1$. 
\end{Remark}
We denote by $E_0 : H^1(S^1; \R^n) \to \R$ the square length functional for
the curves in $\R^n$, namely
$$
E_0(u) = \frac{1}{2} \int_0^1 |\dot{u}(t)|^2 \, dt; \qquad u \in H^1(S^1;
\R^n), 
$$
and more in general, for every $\varepsilon > 0$, we define 
$E_\varepsilon:H^1(S^1,\mathbb{R}^n)$ to be  
$$ 
E_\varepsilon(u) = \int_0^1 \left( \frac{1}{2} |\dot{u}(t)|^2  
- \frac{1}{\varepsilon} V(u(t)) \right) dt, 
\qquad  \forall u \in H^1(S^1;\R^n). 
$$  
The critical points of $E_\e$ are precisely the solutions of problem
\eqref{eq:pe}.

Now we introduce some final notations. Given a covariant tensor $T$ and a
vector field $X$, we denote by  $\mathcal{L}_X T$ the Lie derivative. If $h \in
H^1(S^1; M)$, we define the functions $Q_h, P_h, B_h : S^1 \to \R$ in the
following way 
$$   
Q_h = |H(h(\cdot)) [\dot{h}(\cdot)]|^2; \qquad 
P_h = \dot{h}(\cdot) \cdot H(h(\cdot)) [\dot{h}(\cdot)]; \qquad 
B_h = b(h(\cdot)), 
$$ 
where we have set, for brevity 
$$
b(x) = V''(x)[n_x,n_x], \qquad  x \in M.
$$
Since $h \in H^1(S^1; M)$ it is easy to check that $Q_h, P_h \in L^1(S^1)$
while, since $h$ is continuous form $S^1$ to $M$, $B_h \in C^0(S^1)$.  
Finally we set 
$$
\ov{H} = \sup \left\{ \|H(x)\| : x \in M\right\}.
$$

\section{The case of a non degenerate geodesic}
This Section is devoted to prove Theorem \ref{t:non-deg}. The strategy is the 
following: since problems \eqref{eq:per-T} and \eqref{eq:pe} are
equivalent, we are reduced to find critical points of $E_\e$ for $\e$ small. 
In order to do this, we first find some "pseudo" critical points for $E_\e$, in
Lemma \ref{l:e2}, and we prove the uniform invertibility of $D^2E_\e$ at
these points in Lemma \ref{l:inv-un}. Then, in Proposition \ref{p:contr} we use
the Contraction Mapping Theorem to find "true" critical points of $E_\e$.    

To carry out the first step of our procedure, let us consider the non-degenerate 
geodesic $x_0$ in the statement of Theorem \ref{t:non-deg}. 
$x_0$ induces the smooth map $n_{x_0(\cdot)} : S^1 \to S^{n-1}$.
Hence every  curve $y : S^1 \to \R^n$ can be decomposed into two parts, the
first tangential to  $T_{x_0} M$, and the second normal to $T_{x_0} M$
\begin{equation}\label{eq:ytyn}
y^T = y - (y \cdot n_{x_0}) \, n_{x_0}; \qquad y^N = (y \cdot n_{x_0}) \, n_{x_0}.  
\end{equation}
When $y$ is differentiable we can also decompose the derivatives of $y^T$ and
of $y^N$ into a tangential part and  a normal part. Setting $y_n = y
\cdot n_{x_0}$, there results     
\begin{equation}\label{eq:comp-ypt}
\left( \frac{d}{dt} y^T \right)^T = \nabla_{\dot{x}_0} y^T; \qquad   
\left( \frac{d}{dt} y^T \right)^N = - (y \cdot \dot{n}_{x_0}) \, n_{x_0} = 
- H(x_0)[y,\dot{x}_0] \, n_{x_0}; 
\end{equation}
and 
\begin{equation}\label{eq:comp-ypn}
\left( \frac{d}{dt} y^N \right)^N = \dot{y}_n \, n_{x_0}; \qquad   
\left( \frac{d}{dt} y^N \right)^T = (y \cdot n_{x_0}) \, \dot{n}_{x_0} = 
y_n \, H(x_0) \dot{x}_0. 
\end{equation}
Using equations \eqref{eq:ytyn}, \eqref{eq:comp-ypt} and \eqref{eq:comp-ypn}
one can easily deduce that for some constant $C_0 > 0$ depending only on
$x_0$ there holds  
\begin{equation}\label{eq:equiv-s}
\frac{1}{C_0} \cdot \left( \|z_n\|_{H^1(S^1)} + \|z^T\|_{T_{x_0} H^1(S^1; M)} \right) 
\leq \|z\|_{H^1(S^1; \R^n)} \leq 
C_0 \cdot \left( \|z_n\|_{H^1(S^1)} + \|z^T\|_{T_{x_0} H^1(S^1; M)} \right),
\end{equation}
for all $z \in H^1(S^1; M)$. This means that $H^1(S^1; \R^n) \simeq
T_{x_0} H^1(S^1; M) \oplus H^1(S^1)$ and that the two norms 
$\| \cdot \|_{H^1(S^1; \R^n)}$ and $\| \cdot \|_{T_{x_0} H^1(S^1; M)} + \|
\cdot \|_{H^1(S^1)}$ are equivalent.

Now, roughly, we want to solve the equation $\e \cdot \ddot{x}_\e = - V'(x_\e)$
up to the  first order in $\e$: expanding $V'(x_\e) : S^1 \to \R^n$ as 
$$
V'(x_\e(t)) = \e \cdot \alpha(t) + \e^2 \cdot \beta(t) + O(\e^3),  
$$
we have to find $f$ and $g$ such that 
\begin{equation}\label{eq:deo2}
\ddot{x}_0 (t) = - \alpha(t); \qquad \ddot{f} = - \beta(t).
\end{equation}
In order to solve equation \eqref{eq:deo2}, we first find the explicit
expressions  of $\alpha$ and $\beta$ depending on $f$ and $g$. Since $M$ is
assumed to be non-degenerate, the function $V$ can be written as 
$$
V(x) = \frac{1}{2} b(x) \, d_x^2; \qquad d_x = dist(x,M),
$$ 
where $b : \R^n \to \R$ is a smooth negative function, see condition
\eqref{eq:v-rep}. Because of the factor $\frac{1}{\e}$, to expand
$\frac{1}{\e} V'(x_\e)$ up to the first order in $\e$, we need to take into
account the derivatives of $V$ up to the third order. We fix some point
$x_0(t)$, and  we consider an orthonormal frame $(e_1, \dots, e_n)$ in
$x_0(t)$ such that  $e_1, \dots, e_{n-1}$ form an orthonormal basis for
$T_{x_0(t)} M$, and $e_n$ is orthogonal  to $M$. With simple computations one
can easily check that the only non-zero components of  the second and the
third differential of $V$ at $x_0(t)$ are 
\begin{equation}\label{eq:d2d3v1}  
D^3_{inn} V(x_0(t)) = D_i b(x_0(t)); 
\qquad D^3_{ijn} V(x_0(t)) = b(x_0(t)) \cdot D^2_{ij}d(x_0(t)); \qquad 
i,j = 1, \dots n-1;   
\end{equation} 
\begin{equation}\label{eq:d2d3v2} 
D^2_{nn} V(x_0(t)) = b(x_0(t)),
\qquad D^3_{nnn} V(x_0(t)) = 3 D_n b(x_0(t)).
\end{equation}
The second differential of $d_x$ at $x_0(t)$, see \cite{gt} Appendix, 
is given by 
$$
D^2 d_x [v,w]= \sum_{i,j = 1}^{n-1} H_{ij} \, v_i \, w_j; \qquad v, w \in
T_{x_0(t)} M.  
$$
Here the numbers $H_{ij}$ denote the components of $H(x_0(t))$ with respect to the 
basis $(e_1, \dots, e_{n-1})$ of $T_{x_0(t)} M$. 
In particular from the last formula it follows that 
$$
D^3_{ijn} V(x_0(t)) = b(x_0) \cdot D^2_{ij}d(x_0) = b(x_0) \cdot H_{ij} (x_0), 
\qquad i,j = 1, \dots, n-1. 
$$
Hence by expanding $V'(x_0 + \e \, f + \e^2 \, g)$ in powers of $\e$ we get 
\begin{eqnarray}\label{eq:al}
\alpha(t) \cdot z & = & b(x_0(t)) \, f_n \, z_n, \qquad \forall z \in \R^n; \\
\label{eq:be} \beta(t) \cdot z & = & b(x_0) \, g_n \cdot z_n 
+ \frac{1}{2} \cdot \left[ \sum_{i,j= 1}^{n-1} 
b(x_0) \cdot H_{ij} (f_i \, f_j \, z_n + 
f_i \, f_n \, z_j  + f_j \, f_n \, z_i ) \right.   \\
& + & \left. \sum_{i=1}^{n-1} D_i b (x_0) 
(2 f_i \, f_n \, z_n + f_n^2 \, z_i) + 3 \,
D_n b(x_0) \,  f_n^2 \, z_n \right], \qquad \forall z \in \R^n. \nonumber
\end{eqnarray}
Here $f_i, f_n$, etc. denote the components of the vectors with respect to
the basis  $(e_1, \dots, e_n)$. Taking into account \eqref{eq:al} and
\eqref{eq:be}, the equations in \eqref{eq:deo2} become 
\begin{equation}\label{eq:eqfn}
\ddot{x}_0 = - b(x_0) \, f_n; 
\end{equation}
\begin{eqnarray}\label{eq:eqftgn}
\ddot{f} \, \cdot z & = & - b(x_0) \, g_n \, z_n - \frac{1}{2} \left[ \left( b(x_0) 
\sum_{i,j=1}^{n-1} H_{ij} f_i f_j + 2 \sum_{i=1}^{n-1} 
D_i b(x_0) f_i \, f_n + 3 D_n b(x_0) \,
f_n^2  \right) z_n \right. \nonumber \\
& + & \left. \sum_{i=1}^{n-1} 
\left( 2 \sum_{j=1}^{n-1} H_{ij} \, f_j \, f_n + D_i b(x_0) \, f_n^2 
\right) z_i \right], \quad \forall z \in \R^n.
\end{eqnarray}
Equation \eqref{eq:eqfn} can be solved in $f_n$ with 
\begin{equation}\label{eq:xe} 
f_n(t) = a(t) := \frac{1}{b(x_0(t))} \dot{x}_0 (t) \cdot H(x_0(t))  \dot{x}_0 (t)  
= \frac{1}{b(x_0(t))} H(x_0(t))[\dot{x}_0 (t), \dot{x}_0 (t)]. 
\end{equation} 
In fact, since $x_0$ is a geodesic, it turns out that $\left( \ddot{x}_0
\right)^T = \n_{\dot{x}_0} \dot{x}_0 = 0$; moreover by taking $y^T =
\dot{x}_0$ in formula \eqref{eq:comp-ypt}, we can conclude that 
$$
\left( \ddot{x}_0 \right)^T = 0; \qquad \left( \ddot{x}_0 \right)_n 
= - H(x_0)[\dot{x}_0, \dot{x}_0], 
$$
hence \eqref{eq:xe} follows. As far as \eqref{eq:eqftgn} is concerned, we can
write it in variational form, substituting the expression of $f_n$ according
to \eqref{eq:xe} 
\begin{eqnarray}\label{eq:ftgnv}
\int_{S^1} \dot{f} \, \dot{z} & = & \int_{S^1} b(x_0) \, g_n \, z_n +
\frac{1}{2} \int_{S^1} \left( b(x_0) H(x_0)[f^T, f^T] + 2 a(t) \nabla^T b(x_0)
f^T + 3 D_n b(x_0) \, a^2(t)   \right) z_n \nonumber \\
& + & \int_{S^1} a(t) \, b(x_0(t)) \, H(x_0) [f^T, z^T] + \frac{1}{2}
\int_{S^1}  a^2(t) \nabla^T b(x_0) z^T; \qquad \forall z \in H^1(S^1; \R^n). 
\end{eqnarray}
Here $\nabla^T b$ is the tangential derivative of $b$ on $M$. 

The quantity $\int_{S^1} \dot{f} \, \dot{z}$ can be expressed in a suitable way  
by decomposing $f$ and $z$ into their tangent and normal parts. Using  
equations \eqref{eq:comp-ypt} and \eqref{eq:comp-ypn} corresponding to $f$ and $z$, 
and taking into account \eqref{eq:xe} we have 
\begin{eqnarray*}
\int_{S^1} \dot{f} \, \dot{z} & = & \int_{S^1} \nabla f^T \, \nabla z^T + 
\int_{S^1} H[f^T, \dot{x}_0] \, H[z^T,\dot{x}_0] + \int_{S^1} \dot{a} \,
\dot{z}_n +  \int_{S^1} a \, z_n \, H[\dot{x}_0] \cdot H [\dot{x}_0] \\ 
& + & \int_{S^1} z_n \, H[\dot{x}_0, \nabla f^T] + \int_{S^1} a \,
H[\dot{x}_0, \nabla z^T]  - \int_{S^1} H[\dot{x}_0, f^T] \dot{z}_n -
\int_{S^1} H[\dot{x}_0, z^T] \dot{a}.   
\end{eqnarray*}
From equations \eqref{eq:for-g-f2} and \eqref{eq:de2r-l0} it follows that 
$$
D^2 L_0 (x_0)[z^T, f^T] = \int_{S^1} \nabla f^T \, \nabla z^T + 
\int_{S^1} H[f^T, \dot{x}_0] \, H[z^T,\dot{x}_0] - 
\int_{S^1} H[\dot{x}_0, \dot{x}_0] \, H[z^T, f^T], 
$$
hence there holds 
\begin{eqnarray*}
\int_{S^1} \dot{f} \, \dot{z}  & = & D^2 L_0 (x_0)[z^T, f^T] 
+ \int_{S^1} H[\dot{x}_0, \dot{x}_0] \, H[z^T, f^T] + \int_{S^1} \dot{a} 
\, \dot{z}_n  + \int_{S^1} a \, z_n \, H^2[\dot{x}_0, \dot{x}_0] \\ 
 & + & \int_{S^1} z_n \, H[\dot{x}_0, \nabla f^T] + \int_{S^1} a \,
H[\dot{x}_0, \nabla z^T]  - \int_{S^1} H[\dot{x}_0, f^T] \dot{z}_n -
\int_{S^1} H[\dot{x}_0, z^T] \dot{a}.   \end{eqnarray*}
Taking into account that by the definition of $a(t)$, it is  
$$
\int_{S^1} a(t) \, b(x_0(t)) \, H(x_0) [f^T, f^T] = 
\int_{S^1} H(x_0)[\dot{x}_0, \dot{x}_0] \, H(x_0) [f^T, f^T],
$$
equation \eqref{eq:ftgnv} assumes the form 
\begin{eqnarray}\label{eq:ftgnv2}
D^2 L_0 (x_0)[z^T, f^T] \! & \! + \! & \! \int_{S^1} \dot{a} \, \dot{z}_n
+ \int_{S^1} a \, z_n \, H^2[\dot{x}_0, \dot{x}_0] + 
 \int_{S^1} a \, H[\dot{x}_0, \nabla z^T] - \int_{S^1} H[\dot{x}_0, f^T]
\dot{z}_n \nonumber  \\ \! & \! = \! & \! \frac{1}{2} \int_{S^1} \left( b(x_0)
H(x_0)[f^T, f^T] + 2 a(t) \nabla^T b(x_0) f^T +  3 D_n b(x_0) \, a^2(t) - 2
H[\dot{x}_0, \n f^T] \right) z_n  \\ \! & \! + \! & \! \frac{1}{2}
\int_{S^1}  a^2(t) \nabla^T b(x_0) z^T + \int_{S^1} b(x_0) \, g_n \, z_n  +
\int_{S^1} H[\dot{x}_0, z^T] \dot{a};  \qquad \forall z \in H^1(S^1; \R^n). 
\nonumber
\end{eqnarray}
Now we claim that we can find $f^T$ satisfying 
the following conditions 
\begin{equation}\label{eq:orto}
\left( f^T, \dot{x}_0 \right)_{T_{x_0}H^1(S^1;M)} = 0; 
\end{equation}
\begin{eqnarray}\label{eq:nome}
D^2 L_0 (x_0)[z^T, f^T] & = & - \int_{S^1} a \, H[\dot{x}_0, \nabla z^T] 
+ \frac{1}{2} \int_{S^1} a^2(t) \nabla^T b(x_0) z^T
+ \int_{S^1} H[\dot{x}_0, z^T] \dot{a}; \nonumber \\ & &  \mbox{ for all } z^T
\in  T_{x_0} H^1(S^1; M).  
\end{eqnarray}
In fact, since $D^2L_0(x_0)$ is non degenerate on $(\dot{x}_0)^\perp$, there 
exists $f^T$ satisfying \eqref{eq:orto} and satisfying \eqref{eq:nome} for all
$z^T \in (\dot{x}_0)^\perp$. But since $D^2L_0(x_0)[\dot{x}_0, f^T] = 0$ and
also 
$$
- \int_{S^1} a \, H[\dot{x}_0, \nabla \dot{x}_0] + \frac{1}{2} \int_{S^1} a^2(t)
\nabla^T b(x_0) \dot{x}_0 + \int_{S^1} H[\dot{x}_0, \dot{x}_0] \dot{a} = 0, 
$$
as one can check with simple computations, indeed $f^T$ satisfies equation 
\eqref{eq:nome} for all $z^T \in T_{x_0} H^1(S^1; M)$.
We note that from standard regularity theory for ordinary differential
equations, it turns out that $f^T$ is smooth. Then, choosing $g_n$ such that 
\begin{eqnarray}\label{eq:gn}
b(x_0) \, g_n & = & - \frac{1}{2} \left(
b(x_0) H(x_0)[f^T, f^T] + 2 a(t) \nabla^T b(x_0) f^T +  3 D_n b(x_0) \,
a^2(t) - 2 H[\dot{x}_0, \nabla f^T] \right), \nonumber  \\ 
& - & \ddot a + a \, H^2[\dot{x}_0, \dot{x}_0] +
\frac{d}{dt} \left( H[\dot{x}_0, f^T]\right),  
\end{eqnarray}
with the above choice of $f^T$ equation \eqref{eq:ftgnv2} holds true.
In conclusion, we have solved \eqref{eq:eqfn} and \eqref{eq:eqftgn}, so 
also \eqref{eq:deo2} is satisfied.

We can summarize the above discussion in the following Lemma.

\begin{Lemma}\label{l:e2} 
Let $f_n$ and $f^T$ be given by equations \eqref{eq:xe} and \eqref{eq:orto}
and  \eqref{eq:nome} respectively. Let also $g_n$ satisfy \eqref{eq:gn} and
let $g^T \equiv 0$. Then, setting
\begin{equation}\label{eq:xeps}
x_\e = x_0 + \e \, f + \e^2 \, g,
\end{equation}
there exists $C_1 > 0$ such that for $\e$ suffuciently small 
$$
\|DE_\e(x_\e)\| \leq C_1 \cdot \e^2. 
$$
\end{Lemma} 

\Proof
For all $z \in H^1(S^1; M)$ there  results 
$$
\int_{S^1} \dot{x}_\e \cdot \dot{z} - \frac{1}{\e} \int_{S^1} V'(x_\e) [z] = - 
\int_{S^1} \left( \ddot{x}_\e + \frac{1}{\e} V'(x_\e) \right) \cdot z. 
$$
Furthermore by \eqref{eq:deo2} one has 
$$
\ddot{x}_\e + \frac{1}{\e} V'(x_\e) = \ddot{x}_0 + \e \, \ddot{f} + \e^2 \, 
\ddot{g} + \alpha(t) + \beta(t) + O(\e^2),   
$$
hence it follows that 
$$
\|DE_\e(x_\e)\| \leq C_0 \cdot \e^2 \cdot \|z\|; \qquad \mbox{ for all } z \in
H^1(S^1; M).  
$$
This concludes the proof of the Lemma. 
\QED

The group $S^1$ induces naturally an action on the closed curves given by
$\theta : x(\cdot) \to x(\cdot + \theta)$, $\theta \in S^1$. If $x \in
H^1(S^1; M)$, or $x \in H^1(S^1; \R^n)$ is a non-constant map, then in a small
neighbourhood of $x$, the quotients  $H^1(S^1; M)/S^1, H^1(S^1; \R^n)/S^1$ 
are smooth manifolds. 
The functionals $L_0$ and $E_\e$ are invariant under the action of $S^1$, so 
one can expect to have invertibility of $D^2 L_0$ and $D^2 E_\e$ at some point 
only passing to the quotient spaces. 
In particular, taking into account Definition \ref{d:nd} one has that if 
$h$ is a non trivial geodesic, then $D^2 L_0(h)$ has
zero kernel on $T_h (H^1(S^1; M)/S^1$).  

We want to prove the uniform invertibility of $D^2 E_\e(x_\e)$.  
Hence, in Lemma \ref{l:inv-un} and Proposition \ref{p:contr} below, it will be
understood that we are considering the quotients of $H^1(S^1; M)$ and of 
$H^1(S^1; \R^n)$, without writing it explicitly.

\begin{Lemma}\label{l:inv-un}
Let $x_\e$ be given by formula \eqref{eq:xeps}, where $f$ and $g$ are the
functions in the statement of Lemma \ref{l:e2}. Then there exists $\e_0 > 0$
and  $C_2 > 0$ with the following properties.  For $\e \in (0,\e_0)$ the
operator $D^2 E_\e (x_\e)$ is invertible and  \begin{equation}\label{eq:inv-un}
\left\|(D^2 E_\e (x_\e))^{-1} \right\| \leq C_2; \qquad \e \in (0,\e_0). 
\end{equation}
\end{Lemma}

\Proof 
Given $z, w \in H^1(S^1; \R^n)$ there holds 
\begin{equation}\label{eq:d2exe}
D^2 E_\e (x_\e)[z,w] = \int_{S^1} \dot{z} \, \dot{w} - \frac{1}{\e} 
V''(x_\e)[z,w].  
\end{equation}
Reasoning as in the proof of Lemma \ref{l:e2} we can write the relation 
\begin{eqnarray*}
\int_{S^1} \dot{z} \, \dot{w} & = & \int_{S^1} \nabla z^T \, \nabla w^T + 
\int_{S^1} H[z^T, \dot{x}_0] \, H[w^T,\dot{x}_0] + \int_{S^1} \dot{z}_n \, 
\dot{w}_n + \int_{S^1} z_n \, w_n \, H^2[\dot{x}_0, \dot{x}_0] \\ 
& + & \int_{S^1} w_n \, H[\dot{x}_0, \nabla z^T] + \int_{S^1} z_n \, 
H[\dot{x}_0, \nabla w^T] - \int_{S^1} H[\dot{x}_0, z^T] \dot{w}_n
- \int_{S^1} H[\dot{x}_0, w^T] \dot{z}_n.  
\end{eqnarray*}
Moreover, expanding $V''(x_\e)$ one has 
$$
V''(x_\e) = V''(x_0) + \e \, V'''(x_0)[f] + O(\e^2), 
$$
so, taking into account formulas \eqref{eq:d2d3v1} and \eqref{eq:d2d3v2}
one can check that there  exist $C_3 > 0$ and some smooth functions 
$\ov{C}_i : S^1 \to \R$, $i = 1, \dots, n$, such that   
\begin{eqnarray}\label{eq:v1ove}
\left| \frac{1}{\e} V''(x_\e)[z, w] \right. & - & \left. \frac{1}{\e} b(x_0)
\, z_n \, w_n - f_n \, b(x_0) \, \sum_{i,j=1}^{n-1}  H_{ij} (x_0) \, z_i \, w_j -
\sum_{i=1}^{n-1} \ov{C}_i (z_i \, w_n + w_i \, z_n) \right| \nonumber \\ 
& \leq & C_3 \cdot \e \cdot |z| \cdot |w|.  
\end{eqnarray}
Hence, since $f_n$ satisfies equation \eqref{eq:xe}, there results 
\begin{eqnarray}\label{eq:d2eexe}
D^2 E_\e (x_\e)[z,w] & = & D^2 L_0(x_0) [z,w] + \int_{S^1} \dot{z}_n \, 
\dot{w}_n + \int_{S^1} z_n \, w_n \, H^2[\dot{x}_0, \dot{x}_0] 
- \sum_{i =1}^{n-1} \int_{S^1} \ov{C}_i (z_i \, w_n + w_i \, z_n) \nonumber \\ 
& + & \int_{S^1} w_n \, H[\dot{x}_0, \nabla z^T] + \int_{S^1} z_n \, 
H[\dot{x}_0, \nabla w^T] - \int_{S^1} H[\dot{x}_0, z^T] \dot{w}_n
- \int_{S^1} H[\dot{x}_0, w^T] \dot{z}_n \\ & - & \nonumber 
\frac{1}{\e} \int_{S^1} b(x_0) \, z_n \, w_n + O(\e) \cdot \|z\| \cdot \|w\|.
\end{eqnarray}
Since $x_0$ is a non-degenerate critical point for $L_0$, there exist
subspaces  $W^+, W^- \subseteq T_{x_0} H^1(S^1; M)$ with the following
properties  \begin{description}
\item{(i)} $W^+ \oplus W^- = T_{x_0} H^1(S^1; M)$, \qquad $W^+ \cap W^-$ = 0; 
\item{(ii)} $D^2 L_0 (x_0)$ is positive definite (resp. negative definite) on 
$W^+$ (resp. $W^-$). 
\end{description}
Now, taking into account that $H^1(S^1; \R^n)$ can be decomposed as $H^1(S^1;
\R^n) \simeq T_{x_0} H^1(S^1; M) \oplus H^1(S^1)$, we equip it with the
equivalent scalar product $(\cdot, \cdot)_\sharp$ 
$$
(v,w)_\sharp = (v^T,w^T)_{T_{x_0} H^1(S^1; M)} + (v_n,w_n)_{H^1(S^1)}.
$$
We set also 
$$ 
X^+ = W^+ \oplus H^1(S^1); \qquad X^- = W^- \oplus \{0\},  
$$
and we define 
$$
P^+ : H^1(S^1; \R^n) \to X^+; \qquad \qquad P^- : H^1(S^1; \R^n) \to X^-
$$
to be the orthogonal projections, with respect to $(\cdot, \cdot)_\sharp$
onto the subspaces $X^+$ and $X^-$. 

From equation \eqref{eq:d2eexe} there results
\begin{equation}\label{eq:d2ex-}
D^2 E_\e(x_\e)|_{X^-} = D^2 L_0(x_0)|_{W^-} + o(1). 
\end{equation}
So, since $D^2 E_\e(x_\e)$ is self-adjoint it follows that also 
$$
\left( D^2 E_\e(x_\e) v, w\right)_\sharp = \left(v, D^2 E_\e(x_\e) w
\right)_\sharp = o(1) \, \|v\| \, \|w\|; \quad v \in X^+, w \in X^-. 
$$
This implies that 
\begin{equation}\label{eq:d2ex+}
D^2 E_\e(x_\e)|_{X^+} = P^+ \circ D^2 E_\e(x_\e)|_{X^+} + o(1).  
\end{equation}
Our aim is to prove that there exists $C_4 > 0$ and $\e_1 < \e_0$ such that the 
following properties hold
\begin{description}
\item{$(j)$}
$\left( D^2 E_\e(x_\e) y, y\right)_\sharp \leq - C_4 \|y\|^2; \qquad 
y \in X^-, \e \in (0,\e_1)$; 
\item{$(jj)$}
$\left( D^2 E_\e(x_\e) y, y\right)_\sharp \geq C_4 \|y\|^2; \qquad 
y \in X^+, \e \in (0,\e_1)$. 
\end{description}
There results
$$
\left( D^2 E_\e(x_\e) y, y\right)_\sharp = D^2 E_\e (x_\e) [y, y]; \qquad y
\in X^-,  
$$
so condition $(j)$ follows immediately from \eqref{eq:d2ex-} and $(ii)$. 
On the other hand, for $z \in W^+ \oplus H^1(S^1)$, one can write 
$$
\int_{S^1} H[\dot{x}_0,z^T] \, \dot{z}_n = 
- \int_{S^1} z_n \, \frac{d}{dt} \left( H[\dot{x}_0,z^T] \right).
$$
So, inserting this relation into \eqref{eq:d2eexe} one can easily see that
there exists $\hat{C} > 0$ such that   
\begin{equation}\label{eq:pos-1}
D^2 E_\e (x_\e) [z, z] \geq D^2 L_0 (x_0) [z, z] + \int_{S^1}
\dot{z}_n^2  - \frac{1}{\e} \int_{S^1} b_* \, z_n^2 - \hat{C} \left( \int_{S^1}
z_n^2 \right)^{\frac{1}{2}}  \cdot \|z\|.   
\end{equation}
Here we have set $b_* = \sup_{x \in M} b(x) < 0$. 
Given an arbitrary $\delta > 0$, by the Newton inequality there results 
$$
\left( \int_{S^1} z_n^2 \right)^{\frac{1}{2}} \cdot \|z\| \leq  
\cdot \left( \frac{1}{2} \cdot \delta \cdot \|z\|^2 + \frac{1}{2} \cdot
\frac{1}{\delta}  \cdot \int_{S^1} z_n^2 \right), 
$$
hence, by equation \eqref{eq:equiv-s} it follows that 
\begin{equation}\label{eq:pos-2}
\left( \int_{S^1} z_n^2 \right)^{\frac{1}{2}} \cdot \|z\|
\leq \left( \frac{C_0}{2} \cdot \delta \cdot  \|z^T\|_{T_{x_0}H^1(S^1; M)}^2 +
\frac{C_0}{2} \cdot \delta \cdot \|z_n\|^2_{H^1(S^1)} + \frac{1}{2} \cdot 
\frac{1}{\delta} \cdot \int_{S^1} z_n^2 \right). 
\end{equation}
Hence, since $D^2 L_0 (x_0)$ is positive definite on $W^+$, we can choose 
$\delta$ to be so small that  
$$ 
\delta \cdot \hat{C} \cdot C_0 \leq \min \left\{1, \inf \left\{ D^2
L_0(x_0) [w, w] : w \in W^+, \|w\| = 1\right\} \right\}. 
$$
With this choice of $\delta$, equations \eqref{eq:pos-1} and \eqref{eq:pos-2} 
imply the existence of $\e_2 < \e_1$ for which  
\begin{equation}\label{eq:coee}
D^2 E_\e (x_\e) [z, z] \geq \frac{1}{2} \cdot \left( L_0(x_0) [z^T, z^T] + 
\int_{S^1} \dot{z}_n^2 - \frac{1}{\e} \cdot b_* \cdot \int_{S^1} \, z_n^2
\right), \qquad  \e \in (0, \e_2), z \in X^+. 
\end{equation}
Equation \eqref{eq:coee} together with \eqref{eq:d2ex-} implies $(jj)$ and
concludes the Proof of the Lemma.  
\QED

\begin{Proposition}\label{p:contr}
For $\e$ small, problem \eqref{eq:pe} admits an unique solution $\ov{x}_\e$
which satisfies  
\begin{equation}\label{eq:sol-e2}
\|x_\e - \ov{x}_\e\|_{C^0(S^1)} \leq C_5 \cdot \e^2; \qquad \mbox{ for some }
C_5 > 0.  
\end{equation}
\end{Proposition}
\Proof 
We prove the Proposition by using the Contraction Mapping Theorem. Actually we
want to  find $y \in H^1(S^1; \R^n)$ which satisfies 
$$
D E_\e (x_\e + y) = 0; \qquad \|y\|_{H^1(S^1; \R^n)} \leq C_5 \cdot \e^2. 
$$
We can write 
$$
D E_\e (x_\e + y) = D E_\e (x_\e + y) - D E_\e (x_\e) - D^2 E_\e(x_\e)[y] +  
D E_\e(x_\e ) + D^2 E_\e(x_\e)[y].  
$$
Hence it turns out that 
$$
D E_\e (x_\e + y) = 0 \qquad \Leftrightarrow \qquad 
y = F_\e(y),
$$
where $F_\e : H^1(S^1; \R^n) \to H^1(S^1; \R^n)$ is defined by
$$
F_\e(z) := - \left(D^2 E_\e (x_\e)\right)^{-1} \left[ DE_\e(x_\e) -
\left(D E_\e (x_\e + z) - D E_\e (x_\e) - D^2E_\e (x_\e)[z]\right) \right].
$$
We show that the map $F_\e$ is a contraction in some ball $\ov{B}_\rho =
\left\{ z \in  H^1(S^1; \R^n) : \|z\| \leq \rho \right\}$. In fact, if $z \in
\ov{B}_\rho$, by equation \eqref{eq:inv-un} there results  
\begin{equation}\label{eq:no-fez}
\|F_\e(z)\| \leq C_2 \cdot \left( \left\|DE_\e(x_\e)\right\| + \left\|    
DE_\e (x_\e + z) - D E_\e (x_\e) - D^2E_\e (x_\e)[z] \right\| \right).
\end{equation}
With a straightforward calculation one obtains that for all 
$w \in H^1(S^1; \R^n)$ there holds 
$$
DE_\e (x_\e + z)[w] - D E_\e (x_\e)[w] - D^2E_\e (x_\e)[z,w] = 
\frac{1}{\e} \int_{S^1} \left( V''(x_\e)[z,w] - V'(x_\e + z)[w] +  
V'(x_\e)[w] \right). 
$$
Since $V$ is a smooth function, there exists $C_6 > 0$ such that 
$$
\left| V''(x_\e(t))[z,w] + V'(x_\e(t))[w] - V'(x_\e(t) + z)[w] \right| \leq 
C_6 \cdot \|z\|^2 \cdot \|w\|; \qquad w, z \in \R^n \|z\| \leq 1,
t \in S^1, 
$$ 
so it follows that for $\rho$ sufficiently small 
\begin{equation}\label{eq:rho2}
\left\| DE_\e (x_\e + z) - D E_\e (x_\e) - D^2E_\e (x_\e)[z] \right\| \leq 
\frac{1}{\e} \, C_6 \cdot \|z\|_{H^1(S^1; \R^n)}^2; \qquad 
\|z\|_{H^1(S^1; \R^n)} \leq \rho. 
\end{equation}
Hence, by using equations \eqref{eq:no-fez}, \eqref{eq:inv-un} and
\eqref{eq:rho2}, for $\rho$ sufficienlty small there holds 
\begin{equation}\label{eq:br-br2}
\|F_\e(z)\| \leq C_2 \cdot \left( C_1 \, \e^2 + \frac{1}{\e} \, C_6 
\cdot \rho^2 \right); \qquad \|z\| \leq \rho.  
\end{equation}
Now consider two functions $z, z' \in H^1(S^1; \R^n)$: for all 
$y \in H^1(S^1; \R^n)$ there results 
\begin{eqnarray*}
DE(x_\e \!\! & \!\!\! + \!\!\! & \!\! z)[y] - D^2(x_\e)[z,y] - DE(x_\e +
z')[y] + D^2(x_\e)[z',y] \\  & = & \frac{1}{\e} \int_{S^1} \left(
V''(x_\e(t))[z,y] - V'(x_\e(t) + z)[y] - V''(x_\e(t))[z',y] + V'(x_\e(t) +
z')[y] \right).  
\end{eqnarray*} 
So, taking into account that  
\begin{eqnarray*}
V''(x_\e(t))[z] & - & V'(x_\e(t) + z) - V''(x_\e(t))[z'] + V'(x_\e(t) + z') \\
& = & \int_0^1 \left( V''(x_\e(t) + z + s(z' - z)) - V''(x_\e(t)) \right)[z -
z'] ds, \end{eqnarray*}
there results 
\begin{equation}\label{eq:con}
\left| (F_\e(z) - F_\e(z'), w) \right| \leq \frac{1}{\e} \, C_2 \,  
\sup_{t \in S^1, s \in [0,1]} \|V''(x_\e(t) + z + s(z' - z)) - V''(x_\e(t))\|   
\cdot \|z' - z\|_\infty \cdot \|w\|_\infty. 
\end{equation}
Choosing 
\begin{equation}\label{eq:e2}
\rho = C_7 \cdot \e^2, 
\end{equation}
with $C_7$ sufficiently large, by equations \eqref{eq:br-br2} and
\eqref{eq:con} the map $F_\e$ turns out to be a contraction in $\ov{B}_\rho$.
This concludes the proof.  \QED

\begin{Remark}
By equation \eqref{eq:e2}, it follows that $\|y\| = O(\e^2)$. On the other
hand, the proof of Lemma \ref{l:e2} determines uniquely the normal component
$g_n$ of $g$. In other words, this means that the following condition must be
satisfied 
\begin{equation}\label{eq:yne2}
\|y_n\|_{C^0(S^1)} = o(\e^2). 
\end{equation}
Actually we can prove that \eqref{eq:yne2} holds true. In fact, by 
the proof of Proposition \ref{p:contr}, the fixed point $y$ solves  
\begin{equation}\label{eq:eq-y}
y = \left(D^2 E_\e (x_\e)\right)^{-1} z, 
\end{equation}
with 
$$ 
z = - \left[ DE_\e(x_\e) - \left(D E_\e (x_\e + y) - D E_\e (x_\e) - D^2E_\e
(x_\e)[y]\right) \right].
$$
By using equations \eqref{eq:d2ex-}, \eqref{eq:d2ex+} and \eqref{eq:coee} one
can show that $y_n$ satisfies the inequality 
$$
\|y_n\|_{L^2(S^1)} \leq \ov{C} \cdot \e \cdot \|P^+ z\|_{H^1(S^1;\R^n)} \leq 
\ov{C} \cdot \e \cdot \|z\|_{H^1(S^1;\R^n)}, 
$$ 
for some fixed $\ov{C} > 0$. Since $\|z\|_{H^1(S^1;\R^n)} = O(\e^2)$, see
Proposition \ref{p:contr}, from the Interpolation Inequality (see for example
\cite{b}) it follows that 
$$
\|y_n\|_{C^0(S^1)} \leq \ov{C} \, \|y_n\|_{L^2(S^1)}^{\frac{1}{2}} \cdot 
\|y_n\|_{H^1(S^1)}^{\frac{1}{2}} = o(\e^2).
$$
Hence \eqref{eq:yne2} is proved. 
\end{Remark}

\

\noindent
{\it Proof of Theorem \ref{t:non-deg}} 
We define $u_T$ as $u_T(T \cdot) = \ov{x}_\e(\cdot)$, $T \ , \e^2 = 1$, see
Proposition \ref{p:contr}. Property $(i)$ follows immediately from Lemma
\ref{l:e2} and Proposition   \ref{p:contr}. As far as property $(ii)$ is
concerned, we note that  formulas \eqref{eq:xeps} and \eqref{eq:sol-e2} imply
that $\|\ov{x}_\e - x_0\|_\infty = O(\e)$, hence $\frac{1}{\e}
\|V'(\ov{x}_\e)\|_\infty = O(1)$, uniformly in $\e$. This means that
$\|\ddot{\ov{x}}_\e\|_\infty = O(1)$  uniformly in $\e$. The conclusion
follows from the Ascoli Theorem.   \QED

\section{About some linear ODE's}

The purpose of this section is to perform a preliminary study in order to reduce 
the problem, in Section 5, on the manifold $H^1(S^1; M)$. The arguments are 
elementary, and perhaps our estimates are well known, but for the reader's 
convenience we collect here the proofs.

We start by studying the equation 
\begin{equation}\label{eq:lin}
\begin{cases}
\ddot{v}(t) + \l_0 \cdot v(t) = \s(t), & \mbox{ in } [0,2\pi], \\
v(0) = v(2\pi), \quad \dot{v}(0) = \dot{v}(2\pi), & 
\end{cases}
\end{equation}
where $\l_0 \in \R$ is a fixed constant and $\s(t) \in L^1([0,2\pi])$. 
By the Fredholm alternative Theorem, 
problem (\ref{eq:lin}) admits a unique solution if $\l_0$ is not 
an eigenvalue of the associated homogeneous problem. The eigenvalues are 
precisely the numbers $\{k^2\}, k \in \N$. Since the behaviour of the 
solutions of (\ref{eq:lin}) changes qualitatively when $\l_0$ is positive 
or negative, we distinguish the two cases separately. The former ($\l_0 < 0$)
is related to condition \eqref{eq:v-rep}, namely to the repulsive case. The 
latter ($\l_0 > 0$) is instead related to the attractive case.

\

\noindent
\textbf{Case $\l_0 < 0$}

\

\noindent 
Let $G(t)$ be the Green function for problem (\ref{eq:lin}), namely the 
solution $v(t)$ corresponding to $\s(t) = \delta_0(t)$. 
One can verify with straightforward computations that $G(t)$ is given by 
\begin{equation}
G(t) = \frac{1}{2\sqrt{|\l_0|}\sinh(\pi\sqrt{|\l_0|})} \cosh 
\left(\sqrt{|\l_0|}(t-\pi)\right), \qquad t \in [0,2\pi]. 
\end{equation}
The solution for a general function $\s$ is obtained by convolution, namely 
one has 
\begin{equation}
v(t) = \int_{S^1} G(t-s)\s(s) ds, \qquad t \in [0,2\pi]. 
\end{equation}
In particular the following estimate holds
\begin{equation}\label{eq:bd-f-l1}
\|v\|_{L^\infty} \leq \|G\|_{L^\infty} \cdot \|\s\|_{L^1} = \frac{1}{2\sqrt{|\l_0|}} 
\cdot \|\s\|_{L^1}. 
\end{equation}
Furthermore, if $\s \in L^\infty$ one can deduce 
\begin{equation}\label{eq:bd-f-linf}
\|v\|_{L^\infty} \leq \|G\|_{L^1} \cdot \|\s\|_{L^\infty} = \frac{1}{|\l_0|}
\cdot  \|\s\|_{L^\infty}. 
\end{equation}
The last two estimates hold true if, more in general, the constant $\l_0$ is
substituted by a function bounded above by $\l_0$. This is the content of the
following Lemma. 

\begin{Lemma}\label{l:b-neg}
Let $\l(t)$ be a negative continuous and periodic function on $[0,2\pi]$ such
that  $\l(t) \leq \l_0 < 0$, and let $\s(t) \in L^1([0,2\pi])$. Then the
solution $v(\cdot)$ of problem  
\begin{equation}\label{eq:lin-bt}
\begin{cases}
\ddot{v}(t) + \l(t) v(t) = \s(t), & \mbox{ in } [0,2\pi], \\
v(0) = v(2\pi), \quad \dot{v}(0) = \dot{v}(2\pi). & 
\end{cases}
\end{equation}
satisfies the estimates (\ref{eq:bd-f-l1}) and (\ref{eq:bd-f-linf}). 
\end{Lemma}

\Proof 
The existence (and the uniqueness) of a solution is an easy consequence
of the Lax-Milgram Theorem. 
Let $\ov{v}(\cdot)$ denote the unique solution of (\ref{eq:lin-bt})
corresponding to $\l \equiv \l_0$. We start by supposing that $f \geq 0$. In
this way, by the maximum principle, it must be $v(t), \ov{v}(t) \leq 0$ for all
$t$. Define $y(t) = v(t) - \ov{v}(t)$: it follows immediately by subtraction
that $y$ is a $2 \pi$-periodic solution of the equation  
\begin{equation}\label{eq:diff-sol}
\ddot{y}(t) = \l_0 \ov{v}(t) - \l(t) v(t) \leq |\l_0| \cdot y(t).  
\end{equation}
We claim that it must be $y(t) \geq 0$ for all $t$. Otherwise, there is 
some $t_0$ for which $y(t_0) < 0$, and $y'(t_0) = 0$, since the 
function $y$ is periodic. It then follows from (\ref{eq:diff-sol}) that 
$y(t)$ should be strictly decreasing in $t$, contradicting its periodicity. 
Hence we deduce that 
\begin{equation}\label{eq:bt-b0}
\ov{v}(t) \leq v(t) \leq 0, \qquad \mbox{ for all } t \in [0,2\pi].
\end{equation}
For a general $\s$, we write $\s = \s^+ - \s^-$, where $\s^+$ and $\s^-$ are
respectively the positive and the negative part of $\s$. 
Let also $v^{\pm}, \ov{v}^{\pm}$ denote the solutions corresponding
to $\s^{\pm}$. 
By linearity it is $v(t) = v^+(t) - v^-(t)$ and $\ov{v}(t) = \ov{v}^+(t) -
\ov{v}^-(t)$, so, since $v^{\pm}$ and $\ov{v}^{\pm}$ have definite sign, it
turns out that 
$$
\|v\|_{L^\infty} \leq \max \{ \|v^+\|_{L^\infty} , \|v^-\|_{L^\infty}  \}  \leq \max \{
\|\ov{v}^+\|_{L^\infty} , \|\ov{v}^-\|_{L^\infty}  \} \leq \|\ov{v}\|_{L^\infty} .
$$
This implies immediately the estimates (\ref{eq:bd-c-f-l1}) and 
(\ref{eq:bd-c-f-linf}) for $v(t)$.
\QED

\noindent
We want to prove that the estimates in (\ref{eq:bd-f-l1}) and in 
(\ref{eq:bd-f-linf}) are stable under bounded $L^1$ perturbations of the 
function $\l$. Precisely we consider the following problem, where $\g \in 
L^1(S^1)$. 
\begin{equation}\label{eq:lin-c}
\begin{cases}
\ddot{v}(t) + (\l + \g(t)) \cdot v(t) = \s(t), & \mbox{ in } [0,2\pi], \\
v(0) = v(2\pi), \quad \dot{v}(0) = \dot{v}(2\pi), & 
\end{cases}
\end{equation}
for which it is well known the existence and the uniqueness of a solution $v(t)$. 
\begin{Lemma}\label{lem:stab-f-l1}
Let $A > 0$ be a fixed constant, $\l_0 < 0$ and $\l(t) \in C^0(S^1)$
satisfy $\l(t) \leq \l_0$ for all $t$. Let also $\g \in L^1(S^1)$: 
then, given any number $\d > 0$, if $|\l_0|$ is
sufficiently large the solution $v(\cdot)$ of \eqref{eq:lin-c} satisfies the
inequality   \begin{equation}\label{eq:bd-c-f-l1} 
\|v\|_{\infty} \leq (1 + \delta) \cdot 
\frac{1}{2 \sqrt{|\l_0|}} \cdot \|\s\|_{L^1}.   
\end{equation}  
\end{Lemma}
\Proof
The solution $v$ satisfies the equation
$$
\ddot{v}(t) + \l \, v(t) = \s(t) - \g(t) v(t), \qquad \mbox{ in } [0,2\pi],
$$
with periodic boundary conditions, then by Lemma \ref{l:b-neg} there holds 
$$
\|v\|_{L^\infty} \leq \frac{1}{2\sqrt{|\l_0|}} \cdot \left( 
\|\s\|_{L^1} + M \cdot \|v\|_{L^\infty} \right) \leq \frac{1}{2\sqrt{|\l_0|}} 
\cdot \left( \|\s\|_{L^1} + A \cdot \|v\|_{L^\infty} \right). 
$$ 
This implies immediately (\ref{eq:bd-c-f-l1}) and concludes the proof.
\QED

\begin{Lemma}\label{lem:stab-f-linf}
Let $A > 0$ be a fixed constant, let $\l_0 < 0$, and let $\l(t) \in C^0(S^1)$
satisfy $\l(t) \leq \l_0$ for all $t$. Suppose also that $\s \in
L^\infty(S^1)$. Then given any number $\d > 0$, if $|\l_0|$
is sufficiently large the solution $v(\cdot)$ of problem \eqref{eq:lin-c}
satisfies the inequality  
\begin{equation}\label{eq:bd-c-f-linf}
\|v\|_{L^\infty} \leq (1 + \d) \cdot \frac{1}{|\l_0|} \cdot \|\s\|_{L^\infty}. 
\end{equation} 
\end{Lemma}
\Proof
Let $v_0(t)$ be the solution of problem (\ref{eq:lin-c}) corresponding to $\g
\equiv 0$, so in particular, by Lemma \ref{l:b-neg}, there holds 
 
\begin{equation}\label{eq:x-0} 
\|v_0\|_{L^\infty} \leq \frac{1}{|\l_0|} \cdot \|\s\|_{L^\infty}. 
\end{equation}
Let $y: [0,2\pi] \to \R$ be defined by $y(t) = v_0(t) - v(t)$. By 
subtraction one infers that $y(t)$ is a solution of the problem 
$$
\begin{cases}
\ddot{y}(t) + b y(t) = c(t) v(t), & \mbox{ in } [0,2\pi], \\
y(0) = y(2\pi), \quad \dot{y}(0) = \dot{y}(2\pi). & 
\end{cases}
$$
Hence, by applying inequality (\ref{eq:bd-f-l1}) one deduce sthat 
\begin{equation}\label{eq:inf-y} 
\|y\|_{L^\infty} \leq \frac{1}{2\sqrt{|\l_0|}} \cdot \|\g(\cdot) \cdot v(\cdot) 
\|_{L^1} \leq \frac{1}{2\sqrt{|\l_0|}} \cdot \|v\|_{L^\infty} \cdot \|\g\|_{L^1}.  
\end{equation} 
So, since by Lemma \ref{lem:stab-f-l1} the function  $v(t)$ satisfies 
inequality (\ref{eq:bd-c-f-l1}), it follows that 
$$
\|y\|_{L^\infty} \leq (1 + \d) \cdot \frac{1}{2 \, |\l_0|} \cdot A \cdot
\|\s\|_{L^\infty}.   
$$
Now, taking into account formulas (\ref{eq:bd-c-f-l1}) and (\ref{eq:x-0}) 
it follows that 
\begin{equation}\label{eq:sti-x-f}
\|v\|_{L^\infty} \leq \|v_0\|_{L^\infty} + \|y\|_{L^\infty} \leq  
\frac{1}{|\l_0|} \cdot \|\s\|_{L^\infty} + (1 + \d) \cdot \frac{1}{2 \, |\l_0|}
\cdot A \cdot \|\s\|_{L^\infty} \leq \frac{(1 + \d)}{|\l_0|} \cdot 
\left( 1 + \frac{A}{2} \right).   
\end{equation}
For $|\l_0|$ large this is a better estimate than (\ref{eq:bd-c-f-linf}), and
inserting it in formula (\ref{eq:inf-y}) we obtain 
$$ 
\|y\|_{L^\infty} \leq \frac{1}{2\sqrt{|\l_0|}} \cdot \|v\|_{L^\infty} \cdot
\|\g\|_{L^1}   
\leq \frac{(1 + \d)}{2} \cdot A \cdot \frac{1}{|\l_0|^{\frac{3}{2}}} \cdot
\left( 1 + \frac{A}{2} \right) \cdot  \|\s\|_{L^\infty}.    
$$ 
Using this estimate in (\ref{eq:sti-x-f}) we finally deduce, if $|\l_0|$ is
sufficiently large  
$$
\|v\|_{L^\infty} \leq \frac{1}{|\l_0|} \cdot \|\s\|_{L^\infty} + 
\frac{(1 + \d)}{2} \cdot A \cdot \frac{1}{|\l_0|^{\frac{3}{2}}} \cdot
\left( 1 + \frac{A}{2} \right) \cdot  \|\s\|_{L^\infty} \leq 
(1 + 2 \, \d) \cdot \frac{1}{|\l_0|} \cdot \|\s\|_{L^\infty}.
$$
This concludes the proof. 
\QED

\
 
\noindent
\textbf{Case $\l_0 > 0$}

\

\noindent 
We recall that the estimates of this case will be applied to the study of the 
attractive case. 
We will always take for simplicity $\l_0$ of the form 
$\l_0 = \left(k + \frac{1}{2}\right)^2$. This is in order to assure that 
$\l_0$ is not an eigenvalue of the problem and that the distance of $\l_0$ from the
spectrum is always of order $\sqrt{\l_0}$. 
Let $\ov{G}(t)$ be the Green function for problem (\ref{eq:lin}), namely the 
solution $v(t)$ corresponding to $\s(t) = \delta_0(t)$. 
One easily verifies that $\ov{G}(t)$ is given by 
\begin{equation}
\ov{G}(t) = \frac{1}{2\sqrt{\l_0}} \sin \left(\sqrt{\l_0} t\right), 
\qquad t \in [0,2\pi]. 
\end{equation}
The solution for a general function $\s$ is obtained again by 
convolution, namely one has 
\begin{equation}
v(t) = \int_{S^1} \ov{G}(t-s)\s(s) ds, \qquad t \in [0,2\pi]. 
\end{equation}
In particular the following estimate can be immediately deduced
\begin{equation}\label{eq:ov-bd-f-l1} 
\|v\|_{L^\infty} \leq \|\ov{G}\|_{L^\infty} \cdot \|\s\|_{L^1} =
\frac{1}{2\sqrt{\l_0}}   \cdot \|\s\|_{L^1}.  
\end{equation} 
If moreover $\s \in L^\infty(S^1)$, one can further deduce 
\begin{equation}\label{eq:ov-bd-f-linf}  
\|v\|_{L^\infty} \leq \|\ov{G}\|_{L^1} \cdot \|\s\|_{L^\infty} =
\frac{2}{\sqrt{\l_0}} \cdot \|\s\|_{L^\infty}.   
\end{equation}

\begin{Remark}
We note that, differently from the case of $\l_0 < 0$, the constant is
changed only by a factor $4$ from \eqref{eq:ov-bd-f-l1} to 
\eqref{eq:ov-bd-f-linf}, and is not by a power of $\l_0$, as in  the
preceding case. However, if $\s \in L^\infty$, \eqref{eq:ov-bd-f-linf} could
be a better estimate than \eqref{eq:ov-bd-f-l1}, since $\|\s\|_{L^1} \leq 2 \pi
\, \|\s\|_{L^\infty}$.      
\end{Remark}

The following analogue of Lemmas \ref{lem:stab-f-l1} and \ref{lem:stab-f-linf}
holds, the proof follows the same arguments. 
 
\begin{Lemma}\label{lem:ov-stab-f-l1} 
Let $A > 0$ be a fixed constant, and suppose $\l_0$ is a constant of the form $\l_0
= (m + 1/2)^2$. Suppose that $\s \in L^1(S^1)$. Then if $\l_0$ is sufficiently
large, problem  
\begin{equation}\label{eq:lin-c2} 
\begin{cases} 
\ddot{v}(t) + (\l_0 + \g(t)) \cdot v(t) = \s(t), & \mbox{ in } [0,2\pi], \\ 
v(0) = v(2\pi), \quad \dot{v}(0) = \dot{v}(2\pi). &  
\end{cases} 
\end{equation} 
possesses an unique solution for all $\g(t) \in L^1(S^1)$ with  $\|\g\|_{L^1}
\leq A$. Moreover, given any number $\d > 0$, if $\l_0$ is sufficiently large
then the solution $v(\cdot)$ satisfies the inequality  
\begin{equation}\label{eq:bd-c-f-l12} 
\|v\|_{L^\infty} \leq (1 + \d) \cdot \frac{1}{2 \, \sqrt{\l_0}} \cdot
\|\s\|_{L^1}.   
\end{equation}  
If moreover $\s \in L^\infty(S^1)$, then for $\l_0$ is sufficiently large there
holds   
\begin{equation}\label{eq:bd-c-f-linf2}  
\|v\|_{L^\infty} \leq (1 + \d) \cdot \frac{2}{\sqrt{\l_0}} \cdot \|\s\|_{L^\infty}.
\end{equation}   
\end{Lemma}

\section{Proof of Theorem \ref{t:es-deg}}

The goal of this section is to prove Theorem \ref{t:es-deg}. Two are the main 
ingredients: the first is the reduction on the manifold $H^1(S^1;  M)$,
treated in Subsection 5.1. The second is the study of the reduced functional, 
carried out Subsection 5.2. 

\subsection{The reduction on $H^1(S^1,M)$}

In this subsection we perform a Lyapunov-Schmidt reduction of problem 
\eqref{eq:per-T} on the manifold $H^1(S^1,M)$ of the closed $H^1$ loops on $M$. 
A fundamental tool are the estimates of Section 4.

Solutions of problem \eqref{eq:pe}, and hence of problem \eqref{eq:per-T}, can
be found as critical points of the functional $E_\e : H^1(S^1; \R^n) \to \R$.

If $u \in \RN$ is in a sufficiently small neighbourhood of $M$, 
say if $dist(u, M) \leq \rho_0$ for some $\rho_0 > 0$, then are 
uniquely defined $h(u) \in M$ and $v(u) \in \R$ such that 
\begin{equation}\label{eq:u-hv} 
u = h(u) + v(u) \cdot n_{h(u)}; \qquad dist(u,M) \leq \rho_0. 
\end{equation}
It is clear that $h$ and $v$ depend smoothly on $u$. In particular, if 
$u(\cdot) \in H^1(S^1; \R^n)$, and if $dist(u(t); M) \leq \rho_0$ for all 
$t$, then $h(u(\cdot))$ and $v(u(\cdot))$ are of class 
$H^1(S^1; M)$ and $H^1(S^1)$ respectively. In the sequel we will often 
omit the dependence on $u$ of $h(u)$ and $v(u)$. 
Viceversa, given $h \in M$, and $v \in \R, |v| \leq \rho_0$, then the point 
$u \in \R^n, u = h + v \cdot n_h$ depends smoothly on $h$ and $v$.

If $u(\cdot) \in H^1(S^1; \R^n)$, then there holds 
\begin{equation}\label{eq:ud-hv}  
\dot{u} = \dot{h} + \dot{v} \cdot n + v \cdot \dot{n} = 
\dot{h} + \dot{v} \cdot n + v \cdot H(h) [\dot{h}]
= (Id + v \, H(h)) [\dot{h}] + \dot{v} \cdot n; \quad \mbox{ a. e. in } S^1. 
\end{equation} 
From the last expression it follows in particular that  
\begin{equation}\label{eq:u-dot-2}
|\dot{u}|^2 = |\dot{h}|^2 + |\dot{v}|^2 + v^2 |H(h) [\dot{h}]|^2 + 
2 v \dot{h} \cdot H(h) [\dot{h}].
\end{equation}
We define $\ov{E}_0, \ov{E}_\e : H^1(S^1; M) \times H^1(S^1) \to \R$ and
$\ov{V} : M \times (-\rho_0,\rho_0) \to \R$ to be 
$$
\ov{E}_{0}(h,v) = E_{0}(u(h,v)), \qquad \ov{E}_{\e}(h,v) = E_{\e}(u(h,v));
\qquad h \in H^1(S^1; M), v \in H^1(S^1), \|v\|_{L^\infty} \leq \rho_0, 
$$
$$
\ov{V}(h,v) = V(u(h,v)); \quad h \in M, v \in (-\rho_0,\rho_0). 
$$
Hence, by means of formula \eqref{eq:u-dot-2}, the functional $\ov{E}_0$
assumes the expression 
\begin{eqnarray*}
\ov{E}_0(h,v) & = & \frac{1}{2} \int_{S^1} 
\left( |\dot{h}|^2 + |\dot{v}|^2 + v^2 |H(h) [\dot{h}]|^2 + 
2 v \dot{h} \cdot H(h) [\dot{h}] \right) \\ & = & L_0(h) + 
\frac{1}{2} \int_{S^1} \left(|\dot{v}|^2 + v^2 |H(h) [\dot{h}]|^2 + 
2 v \dot{h} \cdot H(h) [\dot{h}] \right). 
\end{eqnarray*}
If we differentiate $\ov{E}_0$ with respect to a variation $k \in T_h
H^1(S^1; M)$ of $h$, we obtain 
\begin{eqnarray}\label{eq:der-l-k}
D_h \ov{E}_0(h,v)[k]  =  DL_0(h)[k] &\!\!\! + & \!\!\! 
\int_{S^1} v^2 \left( H(h) [\dot{h}] \cdot H(k) [\dot{k}] + H(h) [\dot{h}]
\cdot \mathcal{L}_k H(h) [\dot{h}] \right) \\ &\!\!\! + &\!\!\!  \int_{S^1} v
\left( 2 H(h)[\dot{h}, \dot{k}] + \dot{h} \cdot \mathcal{L}_k H(h) [\dot{h}]
\right). \nonumber  \end{eqnarray}
Here $\mathcal{L}_kH$, see Notations, denotes the Lie derivative of $H$ in the
direction $k$. Similarly, if we differentiate $\ov{E}_0$ with respect to
a variation $w \in H^1(S^1)$ of $v$, we have    
\begin{equation}\label{eq:der-l-w}
D_v \ov{E}_0(h,v)[w]  =  \int_{S^1} \dot{v} \dot{w} dt + 
\int_{S^1} v \, w \, |H(h) [\dot{h}] |^2 dt + \int_{S^1} w \, \dot{h} \cdot 
H(h) [\dot{h}] dt.
\end{equation}
From equations (\ref{eq:u-hv}) and (\ref{eq:ud-hv}) it follows that if
$u(\cdot) \in H^1(S^1; \R^n)$, and if $dist(u(t),M) \leq \rho_0$ for all $t
\in S^1$, then the condition $D\ov{E}_\e(h,v) = 0$ is equivalent to the
system  
\begin{equation}\label{eq:sys}
\begin{cases}
D_h \ov{E}_\e(h,v)[k] = 0 & \forall k \in T_h H^1(S^1; M), \\ 
D_v \ov{E}_\e(h,v)[w] = 0 & \forall w \in H^1(S^1).  
\end{cases}
\end{equation}
Hence, in order to find critical points of $\ov{E}_\e$ (and hence of $E_\e$),
we first solve the second equation in \eqref{eq:sys}. Then, denoting by $v(h)$
this solution, we solve in $h$ the equation 
$$
D_h \ov{E}_\e(u(h,v(h)))[k] = 0, \qquad \forall k \in T_h H^1(S^1; M). 
$$
By (\ref{eq:der-l-w}), the equation in $v$ $D_w \ov{E}_\e(u(h,v)) = 0$ 
means that $v : S^1 \to \R$ is solution of the following problem  
\begin{equation}\label{eq:der-v0}
\begin{cases}
\ddot{v} - |H(h) [\dot{h}] |^2 v = \dot{h} \cdot H(h) [\dot{h}] - 
\frac{1}{\e} \frac{\partial \ov{V}}{\partial v}(h,v). & \mbox{ in } 
[0,2\pi], \\ v(0) = v(2\pi), \quad \dot{v}(0) = \dot{v}(2\pi). &  
\end{cases}
\end{equation}

\begin{Proposition}\label{p:red-rep}
Suppose that the potential $V$ is of repulsive type, namely that \eqref{eq:v-rep} 
holds, and let $A$ be a fixed positive constant. 
Then there exists $\e_A > 0$ and $C_A > 0$ such that for all $\e \in (0,\e_A)$ 
and for all $h \in L_0^A$ equation (\ref{eq:der-v0}) admits an unique 
solution $v(h)$ which satisfies 
\begin{equation}\label{eq:st-n-v}
\|v(h)\|_{C^0(S^1)} \leq C_A \cdot \sqrt{\e}. 
\end{equation}
Moreover the application $h \to v(h)$ from $L_0^A$ to $C^0(S^1)$ is of class 
$C^1$ and compact, namely if $(h_m)_m \subseteq L_0^A$, 
then up to a subsequence $v(h_m)$ converge in $C^0(S^1)$.  
\end{Proposition}

\Proof
Equation (\ref{eq:der-v0}) can be written it in the form 
($Q_h$, $P_h$ and $B_h$ are defined in Section 2) 
\begin{equation}\label{eq:der-v0-2}
\ddot{v} + \left( \frac{1}{\e} B_h - Q_h \right) \, v
= P_h + \frac{1}{\e} 
\left( B_h v - \frac{\partial \ov{V}}{\partial v}(h,v)  \right), 
\end{equation}
Since $V$ is of repulsive type and since every curve $h \in 
H^1(S^1, M)$ is continuous, the function $B_h$ is a negative continuous
function of $t$ and is bounded above by the negative constant 
$b_* = \sup_{x \in M} b(x)$.  
From Lemma \ref{lem:stab-f-l1} it follows that the resolutive operator 
$\Sigma_{h,\e} : L^1(S^1) \to C^0(S^1)$ for the $2 \pi$-periodic problem 
$$
\ddot{v} + \left( \frac{1}{\e} B_h -  Q_h  \right) \, v = \s(t),
$$
is well defined, whenever $h \in H^1(S^1; M)$, $\s \in L^1(S^1)$ and 
$\e > 0$ is sufficiently small. 
Moreover, by Lemma \ref{lem:stab-f-l1} (resp. Lemma \ref{lem:stab-f-linf}) it 
follows that, given $\d > 0$, the solution of problem \eqref{eq:der-v0} 
satisfies the following estimate, provided $\e$ is small enough   
\begin{equation}\label{eq:e-t-p}
\|\Sigma_{h,\e} \s\|_{L^\infty} \leq (1 + \d) \cdot \frac{1}{2}
\frac{\sqrt{\e}}{\sqrt{|b_*|}}  \cdot \|\s\|_{L^1(S^1)} \qquad \left(
\mbox{resp. }  \|\Sigma_{h,\e} \s\|_{L^\infty} \leq (1 + \d) \cdot
\frac{\e}{|b_*|} \cdot \|\s\|_{L^\infty(S^1)} \right).
\end{equation}

For $v \in C^0(S^1)$, let
$$ 
\s_v = P_h + \frac{1}{\e}  
\left( B_h v - \frac{\partial \ov{V}}{\partial v}(h,v) \right), 
$$  
and let $\Theta_{h,\e} : v \to \Sigma_{h,\e} \s_v$. 
Then, as observed in the Notations it is 
\begin{equation}\label{eq:split-f}
P_h \in L^1(S^1); \qquad 
\frac{1}{\e} \left( B_h v - \frac{\partial \ov{V}}{\partial v}(h,v) \right)
\in L^\infty(S^1).  
\end{equation}
We show that $\Theta_{h,\e}$ is a contraction in some suitable ball
$\mathcal{B}_\rho(0) = \{ v \in C^0(S^1) : \|v\|_{L^\infty} \leq \rho\}$, 
where $\rho$ will be chosen appropriately later.

If $\|v\|_{L^\infty} \leq \rho$, from (\ref{eq:split-f}), \eqref{eq:e-t-p} and 
from the linearity of equation \eqref{eq:der-v0-2} it follows that 
$$
\|\Sigma_{h,\e} \s_v\|_{L^\infty} \leq (1 + \d) \cdot \frac{1}{2}
\frac{\sqrt{\e}}{\sqrt{|b_*|}} \cdot \int_{S^1} P_h + 
(1 + \d) \cdot \frac{\e}{|b_*|} \cdot \frac{1}{\e} \cdot 
\left\| B_h v - \frac{\partial \ov{V}}{\partial v}(h,v) \right\|_{L^\infty}.
$$
By the definition of $\ov{H}$ (see Section 2) it follows that we can estimate 
$\|P_h\|_{L^1(S^1)}$ in this way
\begin{equation}\label{eq:phl1}
\|P_h\|_{L^1(S^1)} = \int_{S^1} \dot{h} \cdot H(h)[\dot{h}] \leq 2 \, \ov{H}
\cdot L_0(h) \leq 2 \, \ov{H} \cdot A. 
\end{equation}
Moreover, since $B_h = \frac{\partial^2 \ov{V}}{\partial v^2}(h,0)$, it
turns out that  
\begin{equation}\label{eq:v-2-ord}
B_h \, v - \frac{\partial \ov{V}}{\partial v}(h,v) = \psi(h,v) \cdot
v^2;  \qquad |v| \leq \rho_0, 
\end{equation}
where $\psi(h,v)$ is a smooth (and hence bounded) function. 
So by equations \eqref{eq:phl1} and \eqref{eq:v-2-ord} there holds  
\begin{equation}\label{eq:br-br}
\|\Sigma_{h,\e} \s_v\|_{L^\infty} \leq (1 + \d) \frac{\sqrt{\e}}{\sqrt{|b_*|}}
\cdot \ov{H} \cdot A + (1 + \d) \frac{1}{|b_*|} \cdot \|\psi\|_{L^\infty} \cdot 
\rho^2.
\end{equation}
Furthermore, if we consider two functions $v, v' \in C^0(S^1)$, it turns out
that 
$$
\| \s_v - \Sigma_{h,\e} \|_{L^\infty} = \frac{1}{\e} \left\| 
\psi(h,v) \, v^2 - \psi(h,v')\, (v')^2 \right\|_{L^\infty}. 
$$
Writing
$$
\psi(h,v) \, v^2 - \psi(h,v')\, (v')^2 = \psi(h,v) \left( v^2 - (v')^2 \right) + 
(v')^2 \left( a(h,v) - \psi(h,v') \right),
$$
one can deduce that  
$$
\left\| \psi(h,v) \, v^2 - \psi(h,v')\, (v')^2 \right\|_{L^\infty} \leq
\|\psi\|_{L^\infty} \cdot \|v + v'\|_{L^\infty} \cdot \|v - v'\|_{L^\infty} + 
\|v'\|_{L^\infty}^2 \cdot \|D \psi\|_{L^\infty} \cdot \|v - v'\|_{L^\infty}. 
$$
Hence, if $\|v\|_{L^\infty}, \|v'\|_{L^\infty} \leq \rho$, then by
formula \eqref{eq:e-t-p} it follows that 
\begin{equation}\label{eq:con-br}
\| \Sigma_{h,\e} \s_v - \Sigma_{h,\e} \s_{v'}) \|_{L^\infty} \leq (1 + \d)
\cdot \frac{1}{b_*} \cdot \left( 2 \rho \cdot \|\psi\|_{L^\infty} + \rho^2
\cdot \|D \psi\|_{L^\infty} \right)   \cdot  \| v - v' \|_{L^\infty}.   
\end{equation}
In conclusion, if we choose $\rho = C \cdot \sqrt{\e}$, with $C > 0$
sufficiently large, it follows from equations \eqref{eq:br-br} and
\eqref{eq:con-br} that $\Theta_{h,\e}$ maps $\mathcal{B}_\rho(0)$ in itself 
and is a contraction. Hence we obtain the existence and the local uniqueness of
the solution. It is a standard fact that $h \to v(h)$ is of class $C^1$.   

It remains to prove the compactness assertion. But this is an immediate
consequence of equation (\ref{eq:der-v0}) which implies that $\ddot{v}_n$ is
bounded in $W^{2,1}(S^1)$, so the compactness follows. This concludes the proof
of the Proposition.  
\QED

\subsection{Study of the reduced functional}

In this Subsection we study the functional $E_\e$, after the reduction of 
Subsection 5.1. Roughly, for $\e$ small, the reduced functional turns out to 
be a perturbation of $L_0$ and the standard min-max arguments can be 
used to find critical points. 
 
Let $A > 0$, and let $h \in L_0^A$. By Proposition \ref{p:red-rep} there
exists $v(h)$ such that 
$D_v \ov{E}_\e(h,v(h)) = 0$. Define $L_\e, G_\e : L_0^A \to \R$ in
the following way 
$$
L_\e(h) = \ov{E}_\e (h,v(h)), \quad G_\e(h) = L_\e(h) - L_0(h); \qquad 
h \in L_0^A, \, \e \in (0, \e_A).   
$$

\begin{Proposition}\label{p:l-pert}
Let $A$ be a fixed positive constant, and let $\e_A$ be given by Proposition 
\ref{p:red-rep}. Then there exists $\ov{\e}_A$, $0 < \ov{\e}_A \leq \e_A$,
such that $L_\e : L_0^A \to \R$ is of class $C^1$. Moreover, there exists 
$\ov{C}_A > 0$ such that 
\begin{equation}\label{eq:l-pert}
|G_\e(h)| \leq \ov{C}_A \cdot \sqrt{\e}, \quad |D G_\e(h)| \leq \ov{C}_A 
\cdot \sqrt{\e}; \qquad \forall h \in L_0^A, \, \forall \e \in (0,\ov{\e}_A). 
\end{equation}   
\end{Proposition}

\Proof
Setting for brevity $v = v(h)$, equation (\ref{eq:der-v0}) assumes the form 
$$ 
\int_{S^1} \dot{v} \, \dot{w} + \int_{S^1} v\, w\, Q_h +  
\int_{S^1} 2 w\, P_h - \frac{1}{\e} \int_{S^1}  
\frac{\partial \ov{V}}{\partial v}(h,v)\, w = 0, \qquad \forall w \in
H^1(S^1).   
$$  
In particular, by taking $w = v$ in the last formula it turns out that  
$v$ satisfies the relation 
\begin{equation}\label{eq:var-v} 
\int_{S^1} \dot{v}^2 + \int_{S^1} v^2 \, Q_h +  
\int_{S^1} 2 v \, P_h = \frac{1}{\e} \int_{S^1}  
\frac{\partial \ov{V}}{\partial v}(h,v) \, v. 
\end{equation} 
Since $\ov{E}_\e(h,v)$ has the expression  
$$ 
E_\e(u(h,v)) = \frac{1}{2} \left( \int_{S^1} |\dot{h}|^2 + |\dot{v}|^2 + v^2 \,  
Q_h + 2 v \, P_h \right) - \frac{1}{\e} \int_{S^1} \ov{V}(h,v) dt, 
$$ 
we can take into account formula (\ref{eq:var-v}) 
to obtain  
\begin{equation}\label{eq:ge} 
G_\e(h) = \frac{1}{\e} \int_{S^1} \left( \frac{1}{2} v \cdot  
\frac{\partial \ov{V}}{\partial v}(h,v) - \ov{V}(h,v) \right).  
\end{equation} 
Since $\ov{V}$ is smooth, it turns out that 
$$
\frac{1}{2} v \cdot \frac{\partial \ov{V}}{\partial v}(h,v) - \ov{V}(h,v) = 
\ov{\psi}(h,v) \cdot v^3, \qquad |v| \leq \rho_0,
$$
for some regular function $\ov{\psi}$. Hence, from equations \eqref{eq:st-n-v}
and \eqref{eq:ge} one infers that 
$$ 
|G_\e(h)| \leq \frac{1}{\e} \cdot \|\ov{\psi}\|_{L^\infty} \cdot \|v(h)\|_{L^\infty}^3 \leq 
C_A^3 \cdot \|\ov{\psi}\|_{L^\infty} \cdot \sqrt{\e}.  
$$ 
This proves the first inequality in \eqref{eq:l-pert}.

To show that $G_\e$ is of class $C^1$, we start proving that for some 
$\hat{C}_A > 0$ and for $\e$ sufficiently small, 
\begin{equation}\label{eq:reg-sig}
\ov{\Sigma}_{\e} : L_0^A \to C^0(S^1), \quad h \to v_h \quad 
\mbox{ is of class } C^1, \mbox{ and } \|D_h \ov{\Sigma}_{\e} (h) \| 
\leq \hat{C}_A \cdot \sqrt{\e}. 
\end{equation}
Consider two elements $h$ and $l$ of $H^1(S^1; M)$ and the corresponding
solutions $v(h)$ and $v(l)$ (which for brevity we denote with $v_h$ and 
$v_l$) of problem (\ref{eq:der-v0}). 
By subtraction of the equations there results 
\begin{eqnarray}\label{eq:diff-sols}
\ddot{y} & + & \left( \frac{1}{\e} B_h - Q_h \right) \, y +  
\left( \frac{1}{\e} \Delta B + \Delta Q \right) \, v_l = \Delta P \\ \nonumber
& + & \frac{1}{\e} \cdot (v_h + v_l) \cdot \psi(h,v_h) \, y +  
\frac{1}{\e} \cdot  \left[ \Delta_1 \psi + \Delta_2 \psi \right] \, v_l^2, 
\end{eqnarray}
where we have set 
$$
y  =  v_h - v_l;  \qquad \Delta B = B_h - B_l; \qquad \Delta Q =
Q_h - Q_l;  
$$
$$
\Delta P  =  P_h - P_l;
\qquad 
\Delta_1 \psi = \psi(h,v_h) - \psi(l,v_h); \qquad \Delta_2 \psi = 
\psi(l,v_h) - \psi(l,v_l). 
$$
Since $\psi$ is smooth, there holds 
$$
\Delta_2 \psi = \psi(l,v_h) - \psi(l,v_l) = \tilde{\psi}(l, v_h, v_l) \cdot y, 
$$
for another smooth function $\tilde{\psi}$. Hence we can write 
equation \eqref{eq:diff-sols} in the form 
\begin{eqnarray}\label{eq:gigi} 
\ddot{y} & + & \left( \frac{1}{\e} B_h - Q_h 
- \frac{1}{\e} \cdot (v_h + v_l) \cdot \psi(h,v_h) - \tilde{\psi}(l, v_h, v_l) 
\right) \, y \\ & = & - \left( \frac{1}{\e} \Delta B + \Delta Q \right) \, v_l
+ \Delta P +   \frac{1}{\e} \cdot \Delta_1 a \, v_l^2.   \nonumber
\end{eqnarray} 
Now fix $h \in H^1(S^1; M)$, and let $l \to h$; we can suppose that both 
$h, l \in L_0^A$. If $\e \in (0,\e_A)$ (see Proposition \ref{p:red-rep}), 
then $v_h$ and $v_l$ satisfy inequality \eqref{eq:st-n-v}. Hence by choosing
$\e$ sufficiently small there results 
$$
\frac{1}{\e} B_h - \frac{1}{\e} (v_h + v_l) \cdot \psi(h,v_h) - 
\tilde{\psi} (l, v_h, v_l) \leq \frac{1}{2} \frac{b_*}{\e}. 
$$
Moreover, with computations similar to \eqref{eq:phl1} one can easily prove
that $\|Q_h\|_{L^1(S^1)} \leq 2 \ov{H}^2 \cdot A$, so Lemmas
\ref{lem:stab-f-l1} and  \ref{lem:stab-f-linf} yield, for $\d$ fixed and $\e$
sufficiently small  \begin{eqnarray*} 
\|y\|_{L^\infty} & \leq & (1 + \d) \cdot \frac{1}{2}
\frac{\sqrt{2 \e}}{\sqrt{|b_*|}} \left( C_A \cdot \sqrt{\e} \cdot
\|\Delta Q\|_{L^1(S^1)} \right) \\
& + & (1 + \d) \cdot \frac{2 \e}{|b_*|} \left( C_A \cdot 
\sqrt{\e} \cdot \frac{1}{\e} \cdot \|\Delta B\|_{C^0(S^1)} +
C_A^2 \cdot \frac{1}{\e} \cdot \e \cdot \|\Delta_1 \psi\|_{C^0(S^1)}
\right). 
\end{eqnarray*}

The last quantity tends to $0$ as $l \to h$. This proves the continuity 
of $\ov{\Sigma}_{\e} : h \to v(h)$ from $L_0^A$ to $C^0(S^1)$. The 
differentiability follows from the same reasoning, taking into account 
that the maps 
$$
h \to Q_h \quad \mbox{ and } \quad h \to P_h \quad \mbox{ from } H^1(S^1; M) 
\mbox{ to } L^1(S^1; M) 
$$
$$
h \to B_h \quad \mbox{ from } H^1(S^1; M) \mbox{ to } C^0(S^1; M) \quad 
\mbox{ and } \quad \psi(h,v) \mbox{ from } M \times \left( -\rho_0, \rho_0
\right) \mbox{ to } \R $$  
are differentiable. We remark that the differentials of these functions are
uniformly  bounded for $h \in L_0^A$. Passing to the incremental ratio in 
\eqref{eq:gigi}, one obtains that the directional derivative 
$D \ov{\Sigma}_\e(h) [k]$ along the direction $k \in T_{v_h} H^1(S^1; M)$ is 
the unique $2 \pi$-periodic solution of the equation (linear in $y$) 
\begin{eqnarray}\label{eq:d-ed-var}
\ddot{y} & + & \left( \frac{1}{\e} B_h - Q_h - \frac{2}{\e}
\psi(h,v) \, v - \frac{1}{\e} v^2 D_v \psi(h,v) \right) \, y \\  & = &
\nonumber DP_h[k] + \frac{1}{\e} D_h \psi(h,v)[k] \, v^2 - \frac{1}{\e} 
D B_h[k] \, v - DQ_h[k] \, v. 
\end{eqnarray}
From this formula (using for example a local chart on $H^1(S^1; M)$) one
deduces immediately that the differential $D \ov{\Sigma}_\e$ is continuous on
$L_0^A$.  Now we prove the estimate in \eqref{eq:reg-sig}. 
From equation \eqref{eq:st-n-v} and from the uniform boundedness on the
differentials of $P_h$ and $Q_h$ one can deduce that, for some constant $C$
independent on $h \in L_0^A$ there holds 
$$
\|D P_h [k]\|_{L^1(S^1)} \leq C \cdot \|k\|; \qquad 
\left\| \frac{1}{\e} \cdot v^2 \cdot D_h \psi(h,v) [k]
\right\|_{C^0(S^1)} \leq  C \cdot \|k\|; 
$$
$$
\left\| \frac{1}{\e} \cdot v \cdot D B_h[k] \right\|_{C^0(S^1)} \leq C \cdot
\frac{1}{\sqrt{\e}} \cdot \|k\|; \qquad  
\left\| v \, D Q_h[k] \right\|_{L^1(S^1)} \leq C \cdot \sqrt{\e} \cdot\|k\|.  
$$
Moreover the coefficient of $x$ in formula (\ref{eq:d-ed-var}) for $\e$
sufficiently small can be estimated as 
$$
\frac{1}{\e} B_h - Q_h - \frac{2}{\e} \psi(h,v) \, v - \frac{1}{\e} v^2 D_v
\psi(h,v) \leq \frac{1}{2} \cdot \frac{b_*}{\e},
$$
Hence one can apply Lemmas \ref{lem:stab-f-l1} and \ref{lem:stab-f-linf} to 
obtain 
$$
D \ov{\Sigma}_\e(h) [k] \leq \hat{C}_A \cdot \sqrt{\e} \cdot \|k\|, \qquad 
\forall h \in L_0^A, \forall k \in T_h H^1(S^1; M),
$$
which is the estimate in \eqref{eq:reg-sig}. 

To prove the second inequality in \eqref{eq:l-pert} one can 
write $G_\e$ in the form  $G_\e (h) = \frac{1}{\e} \int_{S^1} \ov{\psi}(h,v_h)
\cdot v_h^3$,  and differentiate with respect to $h$ 
\begin{equation}\label{eq:d-g-e}
DG_\e(h)[k] = \frac{1}{\e} \int_{S^1} \left( 3 \ov{\psi}(h,v) \, v^2\, Dv(h)[k]
+ v^3 \, \frac{\partial \ov{\psi}}{\partial h}[k] + v^3 \,
\frac{\partial \ov{\psi}}{\partial v} Dv(h)[k] \right). 
\end{equation}  
Now it is sufficient to apply the second inequality in \eqref{eq:reg-sig}. 
This concludes the proof.  
\QED

\begin{Lemma}\label{l:psee}
Let $(u_m)_m \subseteq H^1(S^1; \R^n)$ be a Palais Smale sequence for $E_\e$
which for some $R > 0$ satisfies the condition  
\begin{equation}\label{eq:uninbr}
\|u_m(t)\| \leq R, \qquad \forall \, m \in \N, \forall \,t \in S^1.  
\end{equation}
Then, passing to a subsequence, $u_m$ converges strongly in
$H^1(S^1; \R^n)$.  
\end{Lemma}
\Proof 
To prove this claim we first note that by condition \eqref{eq:uninbr}, the sequence 
of numbers $\int_{S^1} V(u_m)$ is bounded. Moreover there holds 
$$
\frac{1}{2} \int_{S^1} \dot{u}_m^2= E_\e(u_m) + \frac{1}{\e} \int_{S^1}
V(u_m),   
$$
so from the convergence of $E_\e(u_m)$ and condition \eqref{eq:uninbr} one
deduces that $(u_m)_m$ is bounded in $H^1(S^1; \R^n)$. Hence, passing to a
subsequence, $u_m \rightharpoonup u_0$ weakly in $H^1(S^1; \R^n)$, and
strongly in $C^0(S^1; \R^n)$. 

As a consequence one has 
\begin{equation}\label{eq:conv-v}
V(u_m(\cdot)) \to  V(u_0(\cdot)) \quad \mbox{ and } V'(u_m(\cdot)) \to 
V'(u_0(\cdot)) \qquad  \mbox{ uniformly on } S^1. 
\end{equation}  
From the fact $D E_\e(u_m) \to 0$, it follows that 
$$
\int_{S^1} \dot{u}_m \, \dot{z} - \frac{1}{\e} \int_{S^1} V'(u_m)[z] =
o(1) \cdot \|z\|, \qquad  \forall z \in H^1(S^1; \R^n),  
$$
where $o(1) \to 0$ as $m \to +\infty$. In particular, taking as test
function $z = u_m$, one has  
\begin{equation}\label{eq:psun}
\int_{S^1} |\dot{u}_m|^2 - \frac{1}{\e} \int_{S^1} V'(u_m)[u_m] = o(1) \cdot 
\|u_m\| = o(1),  \qquad  m \to + \infty.  
\end{equation}
On the other hand, taking $z = u_0$
\begin{equation}\label{eq:psu0}
\int_{S^1} |\dot{u}_0|^2 - \frac{1}{\e} \int_{S^1} V'(u_0)[u_0] = \lim_m 
\int_{S^1} \dot{u}_m \, \dot{u}_0  - \frac{1}{\e} \int_{S^1} V'(u_m)[u_m] =
0, \qquad  m \to + \infty.  
\end{equation}
From equations \eqref{eq:psun} and \eqref{eq:psu0} it follows that 
$$
\lim_m \int_{S^1} |\dot{u}_m|^2 = \int_{S^1} |\dot{u}_0|^2, 
$$
and hence $u_m \to u_0$ strongly in $H^1(S^1; \R^n)$. 
\QED

\begin{Corollary}\label{c:psle}
The functional $L_\e : L_0^A \to \R$, $\e \in (0,\e_A)$, satisfies the Palais 
Smale condition. 
\end{Corollary}

\Proof
Let $(h_m)_m \subseteq L_0^A$ be a Palais Smale sequence for $L_\e$, namely
a sequence  which satisfies 
\begin{description}
\item{(i)} $L_\e (h_m) \to c \in \R$, \qquad as $n \to + \infty$;
\item{(ii)} $\|D L_\e (h_m)\| \to 0$, \qquad as $n \to + \infty$.
\end{description}
We want to prove that $h_m$ converges up to a subsequence to some function 
$h_0 \in H^1(S^1; M)$. By Lemma \ref{l:psee} it is sufficient to show that 
the sequence $u_m = u(h_m, v_{h_m})$ is a Palais Smale sequence for $E_\e$. 
In fact, by the continuity of $h(u)$ the convergence of $u_m$ implies that of
$h_m$. 

To prove that $u_m$ is a Palais Smale sequence for $E_\e$, it is sufficient
to take into account that, by the choice of $v_{h_n}$, it is 
$D_v \ov{E}_\e (h_n,v_{h_n}) = 0$, so one has    
$$
D_h L_\e(h)[k] = D_h \ov{E}_\e (h,v_h)[k] + D_v \ov{E}_\e (h,v_h)[Dv(h)[k]] = 
D_h \ov{E}_\e (h,v_h)[k]. 
$$
Hence, it turns out that 
$$
D E_\e (u_n)[w] = D_h \ov{E}_\e (h_n,v_{h_n})[(Dh(u_n)[w]] + 
D_v \ov{E}_\e (h_n,v_{h_n}) [D v(u_n)[w])] = D L_\e (h_n) [Dh(u_n)[w]].    
$$
It is immediate to check that the differential $D h(u)$ is uniformly bounded
for $|v| \leq \rho_0$, so $\|D L_\e(h_n)\| \to 0$ implies that also
$\|D E_\e(u_n)\| \to 0$. This conclues the proof.  \QED

Since $L_\e$ is of class $C^1$, then it is possible to prove, see for example 
\cite{c}, that there exists a pseudo gradiet $\Omega_\e$ for $L_\e$, namely a 
$C^{1,1}$ vector field which satisfies the conditions 
$$
(\Omega_\e(h),-DL_\e(h)) \geq \|DL_\e(h)\|^2, \qquad 
\|\Omega_\e\| \leq 2 \, \|DL_\e\|; \qquad \forall \, h \in L_0^A.  
$$ 
$\Omega_\e$ induces locally a flow $\phi^t_\e$, $t \geq 0$, for which 
$L_\e$ is non descreasing, and which is strictly decreasing whenever 
$DL_\e \neq 0$.

\begin{Lemma}\label{l:gf}
Let $A > 0$ be a fixed constant. Then there exists $\tilde{\e}_A$ such that 
for all $\e \in (0, \tilde{\e}_A)$ and for all $h \in L_0^A$ the flow
$\phi_t^\e \, h$ is defined for all $t \geq 0$. 
\end{Lemma}  

\Proof 
By Proposition \ref{p:red-rep} there exists $\e_{2A}$ such that for 
$\e \in (0,\e_{2A})$ the functional $L_\e$ is well defined on $L_0^{2A}$. 
Moreover, by Proposition \ref{p:l-pert} we can choose $\ov{\e}_{2A} \leq 
\e_{2A}$ such that 
\begin{equation}\label{eq:bd-G}
|G_\e(h)| \leq \ov{C}_{2A} \cdot \sqrt{\e},  \qquad  
\|DG_\e\| \leq \ov{C}_{2A} \cdot \sqrt{\e}; \qquad \e \in (0,\ov{\e}_{2A}), \, 
h \in L_0^{2A}.  
\end{equation}
We choose $\tilde{\e}_A$ satisfying 
$$
\tilde{\e}_A \leq \ov{\e}_{2A}; \qquad \ov{C}_{2A} \cdot 
\sqrt{\tilde{\e}_A} \leq \frac{1}{4} \, A. 
$$
So, fixed $h \in L_0^{2A}$, by our choice of $\tilde{\e}_A$ and by equation 
\eqref{eq:bd-G}, there holds  
$$
L_0 (h) \leq A \quad \Rightarrow \quad L_\e (h) \leq A 
+ \frac{1}{4} \, A \quad \Rightarrow \quad L_\e (\phi_t^\e \, h) 
\leq A + \frac{1}{4} \, A \quad \Rightarrow \quad L_0 
(\phi_t^\e \, h) \leq A + \frac{1}{2} \, A < 2\, A,
$$
for $\e \in (0, \tilde{\e}_A)$ and whenever $\phi_t^\e \, h$ is well defined. 
Hence $\phi_t^\e \, h$ belongs to the domain of $L_\e$  
whenever $\phi_t^\e \, h$ is well defined. 

By the H\"older inequality and by equation \eqref{eq:bd-G} there results 
$$
\|D L_\e (h)\| \leq \|D G_\e (h)\| + \|D L_0 (h)\| \leq 
\ov{C}_{2A} \cdot \sqrt{\tilde{\e}_A} + \int_{S^1} |\dot{h}|^2 \leq 
\ov{C}_{2A} \cdot \sqrt{\tilde{\e}_A} + A. 
$$
Hence the curve $t \to \phi_t^\e \, h$ is globally Lipshitz, so by the local 
existence and uniqueness of the solutions it can be extended for all 
$t \geq 0$. 
\QED

Now we are in position to prove Theorem \ref{t:es-deg}; the arguments rely on
classical topological methods, see \cite{k}, which we recall for the reader's
convenience. 

\

\noindent
{\em Proof of Theorem \ref{t:es-deg}}. 
Suppose first that $\pi_1 (M) \neq 0$, and let 
$[\a]$ be a non trivial element of $\pi_1 (M)$. Define 
$$
\a_0 = \inf_{h \in [\a] \cap H^1(S^1; M)} L_0(h); 
$$
since the curves belonging to $[\a]$ are non contractible, then it is 
$\a_0 > 0$. By Proposition \ref{p:red-rep} there exists $\e_{2\a_0} > 0$ such 
that for $\e \in (0,\e_{2\a_0})$ the functional $L_\e$ is defined on 
$L_0^{2\a_0}$. Hence for $\e$ sufficiently small it makes sense to define also 
$$
\a_\e = \inf_{h \in [\a] \cap H^1(S^1; M)} L_\e(h).   
$$
By using standard arguments, based on Corollary \ref{c:psle} and 
Lemma \ref{l:gf}, one can prove that the infimum $\a_\e$ is achieved 
by some critical point $h_\e$ of $L_\e$, and that $\a_{\e} \to \a_0$ as 
$\e \to 0$. Now let $\e_k \to 0$: since by Proposition \ref{p:l-pert} it is 
$\|D G_{\e_k}(h_{\e_k})\| \leq O(\sqrt{\e_k})$, then
$h_{\e_k}$ is a Palais Smale sequence for $L_0$. 
Hence, passing to a subsequence, $h_{\e_k}$ must converge to a critical $x_0$ 
point of $L_0$ at level $\a_0$. 
This concludes the proof in the case $\pi_1 (M) \neq 0$. 

Now consider the case of $\pi_1 (M) = 0$. By Theorem \ref{t:hur} and by 
Remark \ref{r:hur}, there exists $\o :S^{l + 1} \to M$ which is not 
contractible. To the function $\o$ one can associate the map 
$$
F_\o : (B^l, \partial B^l) \to (H^1(S^1; M), C_0(S^1; M)),
$$
where $C_0(S^1; M)$ denotes the class of constant maps from $S^1$ to $M$. 
$F_\o$ is defined in the following way. First, identify the closed $l$-ball 
$B^l$ with the half equator on $S^{l + 1} \subseteq \R^{l + 2}$ given by 
$$
\{y = (y_0, \dots, y_{l + 1}) \in S^{l + 1} : y_0 \geq 0, y_1 = 0 \}.
$$
Denote by $c_p(t)$, $0 \leq t \leq 1$, the circle which starts out from 
$p \in B^l$ orthogonally to the hyper plane $\{ y_1 = 0 \}$ and enters 
the half sphere $\{ y_1 \geq 0 \}$. With this we put 
$$
F_\o = \{ f \circ c_p(t) | 0 \leq t \leq 1 \}. 
$$
The correspondence $\o \to F_\o$ is clearly bijective, and hence one can define
the value
$$
\b_0 = \inf_{F_\iota, \iota \in [\o]} \sup_{p \in S^l} L_0 (F_\iota(p)). 
$$
Since $\o$ is non contractible, it is possible to prove that $\b_0 > 0$ so, as 
above, for $\e$ sufficiently small one can define the quantity 
$$
\b_\e = \inf_{F_\iota, \iota \in [\o]} \sup_{p \in S^l} L_\e (F_\iota(p)). 
$$
The number $\b_\e$ turns out to be a critical value for $L_\e$ and, again, if 
$\e_k \to 0$ we can find critical points $h_{\e_k}$ of $L_{\e_k}$ at level 
$\b_{\e_k}$ converging to some geodesic $x_0$ with $L_0 (x_0) = \b_0$. This
concludes the Proof of Theorem \ref{t:es-deg}.   
\QED

\begin{Remark}\label{r:cl}
In \cite{ci} is studied the system \eqref{eq:per-T} on $\R^2$. It is assumed,
roughly, that $V$ possesses a non-contractible set of maxima, and are proved
some existence and multiplicity results for large-period orbits. 
Dealing with two-dimensional systems, we suppose that $V$ is non-degenerate in
the sense of \eqref{eq:ndv}. In particular, this implies that $M$ is a simple
closed curve. When $M$ is a manifold of maxima for $V$, this corresponds to
the case treated in Theorem \ref{t:non-deg}. The non-degeneracy of $V$ allows
us to describe in a quite precise way the asymptotic beheviour of the
trajectories. Furthermore, according to Remark \ref{r:i}, we could also obtain
existence of an arbitrarily large number of adiabatic limits, since we can
apply the above argument to each different element of $\pi_1 (M)$.  
Note also that we can address the case in which $M$ is a manifold of minima,
see the next Section. 
\end{Remark}

\section{Attractive potentials}

In this section we describe how the arguments of Section 5 can be modified
to handle the attractive case, namey that in which $V''(x)[n_x,n_x]$ is
positive for all $x \in M$. Since attractive potentials may cause some
resonance phenomena, the reduction procedure could fail. To overcome this
difficulty we need to make stronger assumptions on $V$.  In particular we
assume that the following condition holds  \begin{equation}\label{eq:d2-co}
\frac{\partial^2 V}{\partial n_x^2}(x) = b_0 > 0; \qquad  \mbox{ for all } x
\in M.   
\end{equation}
We also set 
\begin{equation}\label{eq:l}
\Lambda = \sup_{x \in M} \left|\frac{\partial^3 V}{\partial n_x^3}(x)\right|.
\end{equation}

\begin{Proposition}
Suppose that the potential $V$ satisfies condition \eqref{eq:d2-co}, and let
$A > 0$ be a fixed constant. Then if $\Lambda$ satisfies 
\begin{equation}\label{eq:dis-ea}
4 \, \Lambda \cdot \ov{H} \cdot A < b_0,
\end{equation}
then there exists $\e_k \to 0$ such that equation (\ref{eq:der-v0}) admits a
solution $v(h)$ for every $m$ sufficiently large and for every $h \in L_0^A$.
Moreover, there exists $C_A > 0$ such that $v(h)$ satisfies 
\begin{equation}\label{eq:st-v-a}
\|v(h)\|_{C^0(S^1)} \leq C_A \cdot \sqrt{\e_k}; \qquad \mbox{ for } m
\mbox{ large } \mbox{ and for } h \in L_0^A. 
\end{equation}
Furthermore the application $h \to v(h)$ from $H^1(S^1; M)$
to $C^0(S^1)$ is compact. 
\end{Proposition}

\begin{Remark}
One can easily verify, using rescaling arguments, that condition
(\ref{eq:dis-ea}) is invariant under translation, dilation and rotation in
$\R^n$ of problem \eqref{eq:per-T}.   
\end{Remark}

\Proof
The proof is very similar to that of Proposition \ref{p:contr} and is again
based on the Contraction Mapping Theorem. We choose $\e_k$ such that  
$$
\frac{b_0}{\e_k} = \left(k + \frac{1}{2} \right)^2.
$$
In this way, we can apply the estimates of Section 4 with $\l_0 = 
\frac{b_0}{\e_k}$; we will use the same notations of Section 5. 

Since now $B_h \equiv b_0$, equation (\ref{eq:der-v0}) can be written as  
\begin{equation}\label{eq:der-v0-3}
\ddot{v} + \left( \frac{b_0}{\e_k} - Q_h \right) \, v = \s_v := 
P_h + \frac{1}{\e_k} \psi(h,v) \cdot v^2. 
\end{equation}
The definition of $\psi$ and \eqref{eq:l} imply 
\begin{equation}\label{eq:psim}
|\psi(x,0)| = \frac{1}{2} \, \left| \frac{\partial^3 \ov{V}}{\partial
n_x^3}(x,0)\right|  \leq \frac{1}{2} \, \Lambda, \qquad \forall \, x \in M. 
\end{equation}
So, given an arbitrary number $\delta > 0$, if one takes $\|v\| \leq
\rho$ with $\rho$ sufficiently small and if $h \in L_0^A$, then by 
equations \eqref{eq:bd-c-f-l12} and \eqref{eq:bd-c-f-linf2} one has
\begin{eqnarray*}\label{eq:pino}
\|\Sigma_{h,\e} \s_v\|_{C^0(S^1)} & \leq & (1 + \d) \cdot
\frac{1}{2}\frac{\sqrt{\e_k}}{\sqrt{b_0}} \cdot \left ( \|P_h\|_{L^1(S^1)}
+ 4 \, \frac{1}{\e_k} \|\psi(h,v) \, v^2\|_{L^\infty} \right)
\nonumber 
\\ & \leq & (1 + \d) \cdot \frac{\sqrt{\e_k}}{\sqrt{b_0}} \cdot \ov{H}
\cdot A +  (1 + \d) \frac{1}{\sqrt{b_0} \, \sqrt{\e_k}} \cdot \Lambda \cdot
\rho^2.  
\end{eqnarray*}
So, choosing $\rho = C \cdot \sqrt{\e_k}$, if the following equation is 
satisfied 
\begin{equation}\label{eq:ris-c}
(1 + \d)\, \frac{1}{\sqrt{b_0}} \cdot \Lambda \cdot C^2 + 
(1 + \d) \, \frac{1}{\sqrt{b_0}} \cdot \ov{H} \cdot A \leq C,
\end{equation}
then $\|\Sigma_{h,\e} \s_v\|_{C^0(S^1)} \leq \rho$, so $\mathcal{B}_\rho$ is 
mapped into itself by $\Theta_{h,\e}$. 
If $\d$ is chosen to be small enough, then \eqref{eq:ris-c} is solvable with 
\begin{equation}\label{eq:est-c}
C = \left( 1 - \sqrt{1 - 4\, \frac{\Lambda \,\ov{H} \, A}{b_0}} \right)
\frac{\sqrt{b_0}}{4 \, \Lambda} + o_\d(1). 
\end{equation}
Here $o_\d(1)$ denotes a quantity which tends to $0$ as $\d$ tends to $0$. 
  
Now we show that $\Theta_{h,\e}$ turns out to be a contraction. In fact, 
given two functions $v, v' \in C^0(S^1)$, there holds 
$$
\s_v - \s_{v'} = \frac{1}{\e} \left( (\psi(h,v) - \psi(h,v')) \, v^2 + 
\psi(h,v') \, (v + v') \, (v - v') \right).  
$$
Hence from formula \eqref{eq:bd-c-f-linf2} it follows that 
\begin{eqnarray*}
\| \Sigma_{h,\e} \s_v - \Sigma_{h,\e} \s_{v'} \|_{C^0(S^1)} & \leq & \left( 
\frac{1}{\e_k} \, \|D \psi\|_{L^\infty} \,C^2 \, \rho^2 \, 2 \, (1 + \d) 
\frac{\sqrt{\e_k}}{\sqrt{b_0}} \right. \\ & + & \left. \frac{2}{\sqrt{b_0}} \, 
\sqrt{\e_k} \, \frac{\Lambda}{\e_k} \left(1 - \sqrt{1 - 8  \, \frac{\Lambda 
\,\ov{H} \, A}{b_0}}\right) \, \frac{\sqrt{b_0}}{2\,A} \, \sqrt{\e_k} + o_\d(1) 
\right) \cdot \| v - v'\|_{C^0(S^1)}. 
\end{eqnarray*}
If $\d$ and $\e_k$ are sufficiently small, then the coefficient of 
$\| v - v'\|_{C^0(S^1)}$ in the last formula is strictly less than $1$, hence 
$F$ is a contraction in $B_\r$. This concludes the proof of the existence.  
The compactness can be proved in the same way as before. 
\QED

About the functional $L_\e(h) = \ov{E}_\e(h,v(h))$ we have the following
analogous of Proposition  \ref{p:l-pert}.

\begin{Proposition}\label{p:l-perta}
Suppose condition \eqref{eq:d2-co} holds true, and suppose that $\Lambda$
satisfies  inequality (\ref{eq:dis-ea}). Then for $\e_k$ sufficiently small the
functional $G_\e$ is of class $C^1$ on $L_0^A$ and there exists 
$\ov{C}_A > 0$ such that 
\begin{equation}\label{eq:l-perta}
|G_\e(h)| \leq \ov{C}_A \cdot \sqrt{\e}, \quad |D G_\e(h)| \leq \ov{C}_A 
\cdot \sqrt{\e_k}; \qquad \forall h \in L_0^A, \, \forall \e \in
(0,\ov{\e}_A).  \end{equation}   
\end{Proposition}

\Proof 
The proof is analogous to that of Proposition \ref{p:l-pert}. The only
difference is the estimate of $\|D\ov{\Sigma}_{\e_k}(h)[k]\|_{C^0(S^1)}$,
namely the norm of the $2 \pi$ periodic solution of equation  
\begin{eqnarray}
\ddot{y} & + & \left( \frac{1}{\e_k} B_h - Q_h - \frac{2}{\e_k}
\psi(h,v) \, v - \frac{1}{\e_k} v^2 D_v \psi(h,v) \right) \, y \\  & = &
\nonumber DP_h[k] + \frac{1}{\e_k} D_h \psi(h,v)[k] \, v^2 - \frac{1}{\e_k} 
D B_h[k] \, v - DQ_h[k] \, v, 
\end{eqnarray}
which under assumption \eqref{eq:d2-co} takes the form 
\begin{eqnarray}\label{eq:d-ed-vara}
\ddot{y} & + & \left( \frac{1}{\e_k} b_0 - Q_h - \frac{2}{\e_k}
\psi(h,v) \, v - \frac{1}{\e_k} v^2 D_v \psi(h,v) \right) \, y \\  & = &
\nonumber \theta(t) := DP_h[k] + \frac{1}{\e_k} D_h \psi(h,v)[k] \, v^2 -
DQ_h[k] \, v  \in L^1(S^1).  
\end{eqnarray}
The study of this equation requires some modifications of the arguments in
Section 4. The reason is that the coefficient of $y$ is not uniformly close
(in $L^1(S^1)$) to the constant function $\frac{1}{\e_k} b_0$. 

Equation \eqref{eq:d-ed-vara} can be written in the form
$$
\ddot{y} + \frac{1}{\e_k} b_0 \, y = \theta + \left(Q_h +
\frac{2}{\e_k} \psi(h,v) \, v + \frac{1}{\e_k} v^2 D_v \psi(h,v) \right) \, y.
$$
So, applying equations \eqref{eq:bd-c-f-l12} and \eqref{eq:bd-c-f-linf2} and
taking into account of \eqref{eq:psim} one can deduce 
\begin{eqnarray*}
\|y\|_{C^0(S^1)} & \leq & (1 + \d) \, \frac{\sqrt{\e_k}}{2 \, \sqrt{b_0}}
\cdot \|y\|_{C^0(S^1)} \cdot \left( \|Q_h\|_{L^1(S^1)} + 4 \, \Lambda \, C \,
\sqrt{\e_k} \frac{2}{\sqrt{\e_k}} + 4 \, \frac{1}{\e_k} \, C^2 \,  \e_k \cdot
\|D \psi\|_{L^\infty} \right) \\ 
& + & (1 + \d) \, \frac{\sqrt{\e_k}}{2 \, \sqrt{b_0}} \cdot \|g\|_{L^1(S^1)}.  
\end{eqnarray*}
If the constant $C$ is given by formula \eqref{eq:est-c}, then one can show
that for $\d$ and $\e_k$ sufficiently small the coefficient of 
$\|x\|_{C^0(S^1)}$ on the right hand side is stricly less than $1$, so
estimate \eqref{eq:l-perta} holds true. This concludes the proof. 
\QED

Corollary \ref{c:psle} holds without changes also for this case, and having
Proposition \ref{p:l-perta} one can easily prove the analogous of Lemma
\ref{l:gf}. These facts allow to apply the reduction on $H^1(S^1; M)$ below
the level $\frac{b_0}{4 \, \ov{H} \, \Lambda}$. As a consequence we have 
the following result.  

\begin{Theorem}\label{t:at-gen}
Suppose $\pi_1(M) \neq 0$ (resp. $\pi_1(M) = 0$), and let $\a_0$ (resp.
$\b_0$) be the value which appears in the proof of Theorem \ref{t:es-deg}. 
Suppose the following condition is satisfied 
$$
4 \, \a_0 \cdot \ov{H} \cdot \Lambda < b_0 \qquad \left( \mbox{ resp. } 
4 \, \b_0 \cdot \ov{H} \cdot \Lambda < b_0 \right).
$$
Then there exists a sequence $T_k \to +\infty$ and 
there exists a sequence of solutions $(u_k)_k$ to problem \eqref{eq:per-T}
corresponding to $T = T_k$ such that up to subsequence $u_k(T_k \cdot)$
converge in $C^0(S^1,\R^n)$. The adiabatic limit of $u_k(T_k \cdot)$ is a 
non trivial closed geodesic $x_0$ on $M$ at level $\a_0$ (resp. $\b_0$).    
\end{Theorem}

\

\noindent
{\em Proof of Theorem \ref{t:attr-es}.} 
It is an immediate consequence of Theorem \ref{t:at-gen}, since condition 
$(ii)$ implies that $\Lambda = 0$. 
\QED


\begin{thebibliography}{99} 
 
\bibitem{ab} A. Ambrosetti, M. Badiale, {\em Homoclinics: Poincar\'e-Melnikov 
type results via a variational approach}, Ann.  Inst.  Henri.  
Poincar\'e Analyse Non Lin\'eaire 15 (1998), 233-252.  

\bibitem{b} Brezis, H.: {\em Analyse Fonctionelle}, Masson 1983. 

\bibitem{c} Chang, K.C.: {\em Infinite-dimensional Morse Theory and Multiple 
Solution Problems}, Birkhauser, 1993. 

\bibitem{ci} Cingolani, S., Lazzo, M.: {\em Multiple periodic solutions for
autonomous conservative systems}, to appear on Top. Meth. Nonlin. Anal. 

\bibitem{gt} Gilbarg, D., Trudinger, N.: {\em Elliptic Partial Differential 
equations of the second order}, Springer 1983. 


\bibitem{k} Klingenberg, W. : {\em Riemannian Geometry}, Walter de Gruyter,
1982.
  
 
\bibitem{lf} Lusternik, L.A., Fet, A.I.: {\em Variational problems on closed
manifolds}, Dokl. Akad. Nauk SSSR (N.S.) {\bf 81} 17-18 (Russian) (1951). 


\bibitem{m} Manton, N.: \emph{A remark on the scattering of BPS monopoles}, 
Phys. Lett. B \textbf{110} (1982), 54-56. 

\bibitem{m2} Manton, N.: {\em Monopole interactions at long range}, 
Phys. Lett. B \textbf{154} (1985). 

\bibitem{s} Spivak, M.: {\em A comprehensive introduction to differential geometry}, 
Publish or Perish, 1977. 

 
\bibitem{u} Uhlenbeck, K.: \emph{Moduli Spaces and Adiabatic Limits}, Notices 
A.M.S. \textbf{42-1} (1998), 41-42.  

\end{thebibliography}
\end{document}